# Prescribing a fourth order conformal invariant on the standard sphere - Part I: a perturbation result

Zindine Djadli, Andrea Malchiodi and Mohameden Ould Ahmedou

*Dedicated to the memory of Stephen M. Paneitz*

**Abstract:** In this paper we study some fourth order elliptic equation involving the critical Sobolev exponent, related to the prescription of a fourth order conformal invariant on the standard sphere. We use a topological method to prove the existence of at least a solution when the function to be prescribed is close to a constant and a finite dimensional map associated to it has non-zero degree.

§**1. Introduction.**

**1.1. Statement of the problem.**

In 1983, Paneitz [38] discovered a particularly interesting fourth order operator on 4-manifolds, given by

$$P_g^4 \varphi = \Delta_g^2 \varphi - div_g \left( \frac{2}{3} Scal_g g - 2Ric_g \right) d\varphi \quad ,$$

where $Scal_g$ denotes the scalar curvature of $(M^4, g)$ and $Ric_g$ denotes the Ricci curvature of $(M^4, g)$. This operator enjoys many interesting properties (in particular, it is conformally invariant), and can be seen as a natural extension of the well known second order conformal operator $\Delta_g$ on 2-manifolds. This operator gives rise to a fourth order conformal invariant ($|Ric_g|$ denoting the norm of $Ric_g$ with respect to the metric)

$$Q_g = \frac{1}{12} \left( \Delta_g Scal_g + Scal_g^2 - 3 |Ric_g|^2 \right) \quad ;$$

namely, under the conformal change $g' = e^{2w} g$, $Q_{g'}$ and $Q_g$ are related by

$$P_g^4 w + 2 Q_{g'} e^{4w} = 2 Q_g \quad .$$

There is also a natural fourth order Paneitz operator on manifolds of dimension greater than 4, introduced by Branson [10] and defined as follows. Given a smooth compact Riemannian $n$-manifold $(M, g)$, $n \geq 5$, let $P_g^n$ be the operator defined by

$$P_g^n u = \Delta_g^2 u - div_g \left( a_n Scal_g g + b_n Ric_g \right) du + \frac{n-4}{2} Q_g^n u \quad ,$$

where

$$a_n = \frac{(n-2)^2 + 4}{2(n-1)(n-2)} \quad , \qquad b_n = -\frac{4}{n-2} \quad ,$$

$$Q_g^n = \frac{1}{2(n-1)} \Delta_g Scal_g + \frac{n^3 - 4n^2 + 16n - 16}{8(n-1)^2(n-2)^2} Scal_g^2 - \frac{2}{(n-2)^2} |Ric_g|^2 \quad .$$

If $\tilde{g} = \varphi^{\frac{4}{n-4}} g$ is a conformal metric to $g$, then for all $u \in C^\infty(M)$ we have

$$P_g^n(u\varphi) = \varphi^{\frac{n+4}{n-4}} P_{\tilde{g}}^n(u) \quad .$$

In particular

(1.0.1) $$P_g^n(\varphi) = \frac{n-4}{2} Q_{\tilde{g}}^n \varphi^{\frac{n+4}{n-4}} .$$

Many interesting results on the Paneitz operator and related topics have been recently obtained by several authors. Among them, let us mention Branson [9], Branson-Chang-Yang [11], Chang-Yang [20], Chang-Gursky-Yang [15], Chang-Qing-Yang [16] and [17], Gursky [27]. See also the surveys Chang [13] and Chang-Yang [18] for more results on basic properties of the Paneitz operator.

In view of relation (1.0.1), a natural problem to propose is to prescribe the conformal invariant $Q$. Equation (1.0.1), where $\varphi$ has to be considered as unknown, has a variational structure, and hence solutions occur as critical points of the associated Euler functional. A natural space to look in for solutions is the Sobolev space $H_2^2(M)$ (see section 2). Due to the non-compactness of the injection of $H_2^2(M)$ into $L^{\frac{2n}{n-4}}(M)$, the Euler functional does not satisfy the Palais-Smale condition, which leads to the failure of the standard critical point theory.

In this paper we consider the problem of prescribing $Q$ on the standard sphere $(\mathbb{S}^n, h)$, $n \geq 5$, where the expression of $P_h^n$ is given by

(1.0.2) $$P_h^n u = \Delta_h^2 u + c_n \Delta_h u + d_n u ,$$

where

(1.0.3) $$c_n = \frac{1}{2}(n^2 - 2n - 4) , \quad d_n = \frac{n-4}{16} n(n^2 - 4) .$$

More precisely, given a smooth positive function $f$ on $\mathbb{S}^n$, can one find a metric $g$ conformal to $h$ for which the $Q_g$-curvature is given by $f$? Setting $g = u^{\frac{4}{n-4}} h$, this amounts to solving the following equation

(E) $$\begin{cases} \Delta_h^2 u + c_n \Delta_h u + d_n u = \frac{n-4}{2} f u^{\frac{n+4}{n-4}} , \\ u > 0 \text{ on } \mathbb{S}^n . \end{cases}$$

We recall that, by the regularity results of Djadli-Hebey-Ledoux [21] (see also Van der Vorst [40]), a weak solution of (E) is indeed a smooth solution.

Let us observe that, in the case of the sphere (in view of the uniqueness result of Lin [35]), a solution of (E) cannot be achieved as a minimum for the extremal problem

$$\inf \left\{ \frac{\fint_{\mathbb{S}^n} (\Delta_h u)^2 + c_n \fint_{\mathbb{S}^n} |\nabla u|^2 + d_n \fint_{\mathbb{S}^n} u^2}{\left( \fint_{\mathbb{S}^n} f |u|^{2^\#} \right)^{\frac{2}{2^\#}}} \; \Big| \; u \in H_2^2(\mathbb{S}^n) \; u \not\equiv 0 \right\} .$$

Moreover, due to the invariance of equation (E) under the group of conformal transformations of $(\mathbb{S}^n, h)$, Kazdan-Warner type obstructions occur. Precisely, the following Theorem holds

**Theorem 1.1.** [21] *Let $(\mathbb{S}^n, h)$ be the standard $n$-dimensional unit sphere, $n \geq 5$, and $f$ a smooth function on $\mathbb{S}^n$. If $u$ is a smooth solution of the equation*

$$P_h^n u = \frac{n-4}{2} f |u|^{\frac{8}{n-4}} u ,$$

*then for any eigenfunction $\xi$ of $\Delta_h$ associated to the first non zero eigenvalue $\lambda_1 = n$ it is*

$$\int_{\mathbb{S}^n} <\nabla f, \nabla \xi> |u|^{2^\#} dv_h = 0 ,$$

where $2^{\#} = \frac{2n}{n-4}$. In particular, for any eigenfunction $\xi$ of $\Delta_h$ associated to the first non zero eigenvalue and for any $\varepsilon \in \mathbb{R}$, equation $(E)$ with $f = 1 + \varepsilon\xi$ admits no positive solution.

In fact, this Theorem is proved in Djadli-Hebey-Ledoux [21] for positive solution, but, with the same proof, one can prove that it remains true for changing sign solutions.

The problem of finding conditions on $f$ such that $(E)$ possesses a solution can be seen as the generalization to the Paneitz operator of the so-called "Nirenberg problem"; namely: which functions $K \in C^\infty(\mathbb{S}^n)$ are the scalar curvature of a metric conformal to the standard one? Theorem 1.1 is the analogue of the classical obstruction for the Nirenberg problem (see Kazdan-Warner [31]). The Nirenberg problem has been studied by several authors; among others, we mention Ambrosetti-Garcia Azorezo-Peral [2], Ambrosetti-Malchiodi [3], Aubin [4], Bahri-Coron [6], Chang-Yang [19], Chang-Gursky-Yang [14], Hebey [28], Li [33], Lin [34], Schoen-Zhang [39].

**1.2. Statement of the results.**

It is well known (see for instance Berger-Gauduchon-Mazet [8]), that denoting

$$\mathbb{H} = \{\xi \in C^\infty(\mathbb{S}^n) \quad \xi \neq 0 \quad | \quad \Delta_h \xi = n\xi\} \quad,$$

it is $\dim \mathbb{H} = n+1$, and $\xi \in \mathbb{H}$ if and only if $\xi$ is the restriction to $\mathbb{S}^n$ of a linear form on $\mathbb{R}^{n+1}$. It follows that, denoting by $(x_i)_{i \in \{1,\ldots,n+1\}}$ the canonical coordinates of $\mathbb{R}^{n+1}$ and by $\xi_i$ the restriction to $\mathbb{S}^n$ of the $i$-th canonical projection associated to these coordinates, $(\xi_i)_{i \in \{1,\ldots,n+1\}}$ is a basis of $\mathbb{H}$. In the sequel of the paper we set

$$\overrightarrow{\xi} = (\xi_1, \ldots, \xi_{n+1}) \in \mathbb{H}^{n+1} \quad.$$

Let us now introduce some notation (following those introduced by Chang-Yang [19]). First we will denote by

$$\fint = \frac{1}{Vol_h(\mathbb{S}^n)} \int \quad.$$

For $P \in \mathbb{S}^n$ and $t \in [1; +\infty)$, we denote by $p = \frac{t-1}{t}P \in B$, $B$ being the unit ball of $\mathbb{R}^{n+1}$. Using stereographic coordinates with projection through the point $P$, we set

$$\varphi_{P,t}(y) = ty \quad.$$

This set of conformal transformations is diffeomorphic to the unit ball $B \subset \mathbb{R}^{n+1}$ with the identity transformation identified with the origin and $\varphi_{P,t} \leftrightarrow \frac{t-1}{t}P = p \in B$. In the sequel of the paper, we will consider the map $G : B \to \mathbb{R}^{n+1}$ given by

$$G(P,t) = \fint_{\mathbb{S}^n} (f \circ \varphi_{P,t}) \overrightarrow{\xi} dv_h \quad,$$

where $\overrightarrow{\xi}$ is defined above (note that, here and in the sequel of the paper, we identify $(P,t) \in \mathbb{S}^n \times [1; +\infty)$ and $p = \frac{t-1}{t}P \in B$).

**Définition 1.2.** Let $f \in C^\infty(\mathbb{S}^n)$. We say that $f$ is uniformly non degenerated of order $\alpha$ if there exists $t_1 > 1$ and $C > 0$ such that for all $t \geq t_1$ and for all $P \in \mathbb{S}^n$

$$|G(P,t)| \geq \frac{C}{t^\alpha} \qquad \text{when } \alpha < n \quad;$$

$$|G(P,t)| \geq \frac{C}{t^n} \log t \qquad \text{when } \alpha = n \quad.$$

Our main Theorem is the following

**Theorem 1.3.** *There exists $\varepsilon(n) > 0$, depending only on the dimension of $\mathbb{S}^n$, such that $(E)$ admits a solution for all functions $f \in C^\infty(\mathbb{S}^n)$ satisfying the following assumptions:*

$$(H1) \quad \left\| f - \frac{n(n^2-4)}{8} \right\|_\infty \leq \varepsilon(n) \quad ;$$

$(H2)$ $\quad f$ *is uniformly non degenerated of order at most $n$ if $n$ is even,*
    *and of order at most $n-1$ if $n$ is odd;*

$(H3)$ $\quad deg\left(G, \{p \mid t < t_0\}, 0\right) \neq 0$ *for some $t_0$ large enough.*

Using Theorem 1.3, one can prove the following Corollary

**Corollary 1.4.** *Assume $f \in C^\infty(\mathbb{S}^n)$, $n \geq 5$, is a Morse function satisfying the condition $(H1)$ of Theorem 1.3 and the following assumptions*

$$|\nabla f(P)| = 0 \implies \Delta_h f(P) \neq 0 \quad ,$$

*and*

$$\sum_{P \in \mathbb{S}^n \mid \nabla f(P)=0 \text{ and } \Delta_h f(P)>0} (-1)^{m(P,f)} \neq -1 \quad ,$$

*where $m(P,f)$ is the Morse index of $f$ at $P$. Then $(E)$ admits a solution.*

As observed in Chang-Gursky-Yang [15], the condition on the index in Corollary 1.4 is stronger than the condition $(H3)$ of Theorem 1.3. Recently, Felli [26] proved an existence result for $(E)$ related to Corollary 1.4 using a method introduced by Ambrosetti and Badiale [1].

### 1.3. Comments.

Theorem 1.3 is a generalization to the Paneitz operator of the results of Chang-Yang [19]. Apart from the difficulties mentionned above, differently from the case of the scalar curvature, there are other specific problems. The first is that obtaining positive solutions is not immediate; in the case of the scalar curvature equation, one can obtain positivity just by using the fact that $|\nabla |u|| = |\nabla u|$ *a.e.*. This is not anymore true for the fourth order equation because of the term $\int_{\mathbb{S}^n} (\Delta_h u)^2$ (and the fact that $|u|$ not necessarily belongs to $H_2^2(\mathbb{S}^n)$). We also encounter the same difficulty in obtaining a priori estimates (see section 5) since we cannot apply the Moser iteration scheme. Our main tool to handle such difficulties is an a priori estimate for fourth order operator (see Proposition 3.1) which goes back to a device introduced by Van der Vorst [40] (see also Lee-Parker [32] and Djadli-Hebey-Ledoux [21]) to prove regularity of solutions. Starting from Proposition 3.1, we propose a way to handle positivity and a priori estimates for the above fourth order operator (see Remark 5.9). Another main ingredient in the proof of Theorem 1.3 is a fourth order Sobolev-Aubin type inequality, stated in Theorem 4.1, which could be interesting also in itself.

In the second part of the paper (Djadli-Malchiodi-Ould Ahmedou [24]), we give some existence results without the "close to a constant condition" in dimension 5 and 6. Part of these results are obtained via a continuity argument whose starting point is Theorem 1.3 (and generalize to the Paneitz operator previous results of Chang-Gursky-Yang [14], Bahri-Coron [6], Li [33] and Benayed-Chen-Chtioui-Hammami [7] for the Nirenberg problem). Our work in the second part of this paper is based on a fine blow-up analysis for a family of equations approximating $(E)$. This analysis extend to this framework the notion of isolated and isolated simple blow-up points introduced by R.Schoen (and used extensively by Y.Y. Li [33]). This blow-up analysis requires some Harnack type inequalities for fourth order operators obtained from the interpretation of higher order equations as systems (from which we extract some useful informations on the sign of the laplacian of the solutions). See Djadli-Malchiodi-Ould Ahmedou [24] for details.

The paper is organized as follows : in section 2 we recall some useful facts, in section 3 we derive some a priori estimates for some fourth order equations and in section 4 we prove the Sobolev-Aubin type inequality (see

Theorem 4.1); in section 5 we perform a finite dimensional reduction of the problem, based on Theorem 4.1. In particular we give an uniform estimate on a family of solutions $u_p$ of this reduced problem, depending on a $(n+1)$-dimensional parameter $p$. In section 6 we show that the map $p \to u_p$ is well defined and continuous. In section 7 we prove that the maps $G$ and $p \to u_p$ are related by some degree formula; in the last section we prove Theorem 1.3 and Corollary 1.4.

**Acknowledgments:** Part of this work was accomplished when the first author enjoyed the hospitality of the SISSA at Trieste. He would like to thank A.Ambrosetti for the invitation and all the members of this institute for their hospitality. It is also a pleasure to thank A.Ambrosetti, A.Chang and P.Yang for valuable discussions and helpful hints. A. M. is supported by a Fulbright Fellowship for the academic year 200-2001 and by M.U.R.S.T. under the national project *"Variational methods and nonlinear differential equations"*. M.O.A. research is supported by a postdoctoral fellowship from SISSA.

## §2. Preliminaries on Paneitz operator and background material.

A natural space when studying the Paneitz operator $P_h^n$ is the Sobolev space $H_2^2(\mathbb{S}^n)$ defined as the completion of $C^\infty(\mathbb{S}^n)$ with respect to the norm

$$\|u\|^2 = \left\|\nabla^2 u\right\|_2^2 + \|\nabla u\|_2^2 + \|u\|_2^2 \quad,$$

where $\|.\|_p$ stands for the $L^p$-norm with respect to the Riemannian measure $dv(h)$. As is well known and easy to see, for all $u \in C^\infty(\mathbb{S}^n)$

$$(\Delta_h u)^2 \le n \left|\nabla^2 u\right|^2 \quad.$$

Conversely, using the Bochner-Lichnerowicz-Weitzenbock formula, there holds

$$\int_{\mathbb{S}^n} \left|\nabla^2 u\right|^2 dv(h) \le \int_{\mathbb{S}^n} |\Delta_h u|^2 dv(h) + k \int_{\mathbb{S}^n} |\nabla u|^2 dv(h) \quad,$$

where $k$ is such that $Ric_g \ge -kg$. Hence $\|.\|_{H_2^2}$ defined by

$$\|u\|_{H_2^2}^2 = \int_{\mathbb{S}^n} (P_h u)\, u\, dv(h)$$

is a norm on $H_2^2(\mathbb{S}^n)$ which is equivalent to the classical one $\|.\|$.

By the Sobolev embedding Theorem, for $n \ge 5$, one gets an embedding of $H_2^2(\mathbb{S}^n)$ in $L^{2^\#}(\mathbb{S}^n)$ where $2^\# = \frac{2n}{n-4}$. This embedding is critical and continuous; so there exists $A \in \mathbb{R}$ such that

$$\|u\|_{2^\#}^2 \le A \int_{\mathbb{S}^n} (P_h u)\, u\, dv(h) \quad.$$

Setting

(2.0.1) $$K_0 = \pi^2 n(n-4)(n^2-4)\Gamma\left(\frac{n}{2}\right)^{\frac{4}{n}} \Gamma(n)^{-\frac{4}{n}} \quad,$$

we have the following Theorem (see Djadli-Hebey-Ledoux [21]) regarding the best constant for the embedding $H_2^2(\mathbb{S}^n) \hookrightarrow L^{2^\#}(\mathbb{S}^n)$.

**Theorem 2.1.** *Let $(\mathbb{S}^n, h)$ be the standard unit sphere of $\mathbb{R}^{n+1}$, $n \ge 5$, and $P_h$ be the Paneitz operator $P_h = \Delta_h^2 + c_n \Delta_h + d_n$. Then*

$$\inf_{u \in C^\infty(\mathbb{S}^n)\setminus\{0\}} \frac{\fint_{\mathbb{S}^n} (P_h u)\, u}{\left(\fint_{\mathbb{S}^n} |u|^{2^\#}\right)^{\frac{2}{2^\#}}} = \frac{1}{K_0} \omega_n^{-\frac{4}{n}} \quad,$$

where $\omega_n$ denotes the volume of the standard unit $n$-sphere.

We also recall the following result of Djadli-Hebey-Ledoux [22]

**Theorem 2.2.** *Let $(M,g)$ be a smooth compact $n$-dimensional Riemannian manifold, $n \geq 2$. Then for any positive real number $\alpha$, the spectrum of the operator $\Delta_g^2 + \alpha \Delta_g$ is exactly $\{\lambda^2 + \alpha\lambda , \lambda \in Sp(\Delta_g)\}$, with the additional property that $u$ is an eigenfunction of $\Delta_g^2 + \alpha\Delta_g$ associated to $\mu = \lambda^2 + \alpha\lambda$ if and only if $u$ is an eigenfunction of $\Delta_g$ associated to $\lambda$.*

## §3. An a priori estimate for a fourth order operator.

In this section, following a previous argument by Van der Vorst [40] (see also Djadli-Hebey-Ledoux [21]), we prove an a priori estimate for solutions some fourth order equation. The result plays the same role as the Moser iteration scheme, which is used to gain integrability in second order equations.

**Proposition 3.1.** *Let $(\mathbb{S}^n, h)$ be the standard $n$-sphere, let $\alpha > \frac{n}{4}$, $p \in [1; 2^\# - 1]$ and let $q \in L^\alpha(\mathbb{S}^n)$, $w \in L^\infty(\mathbb{S}^n)$. Assume that $u \in H_2^2(\mathbb{S}^n)$ is a weak solution of the equation*

$$\left(\Delta_h + \frac{c_n}{2}\right)^2 u = q|u|^{p-1}u + w \quad.$$

*Then there exist a postive constant $\delta^\infty$ depending only on $n$ and $\alpha$, and a positive constant $C$ depending only on $n$, $\alpha$ and $\|q\|_\alpha$ such that if $q|u|^{p-1} \in L^{\frac{n}{4}}(\mathbb{S}^n)$ and $\left\|q|u|^{p-1}\right\|_{\frac{n}{4}} \leq \delta^\infty$ then*

$$\begin{cases} u \in L^\infty(\mathbb{S}^n) \quad; \\ \|u\|_\infty \leq C\left(\|w\|_\infty + \|w\|_\infty^p\right) \end{cases} \quad.$$

**Proof.** Let $L_h = \Delta_h + \frac{c_n}{2}$. For all $s > 1$ and for all $g \in L^s(\mathbb{S}^n)$, there exists an unique function $\varphi \in H_4^s(\mathbb{S}^n)$ such that $L_h^2 \varphi = g$. By standard elliptic arguments and the Sobolev embedding theorem that there exists a constant $C(n,s)$, depending only on $n$ and $s$, such that

(3.1.1) $$\begin{cases} \|\varphi\|_{\frac{ns}{n-4s}} \leq C(n,s) \|g\|_s & \text{if } n > 4s \quad; \\ \|\varphi\|_\infty \leq C(n,s) \|g\|_s & \text{if } n < 4s \quad. \end{cases}$$

Now, let $s > 1$ and let $v \in L^s(\mathbb{S}^n)$. It is clear that $q|u|^{p-1} v \in L^{\hat{s}}(\mathbb{S}^n)$ with $\hat{s} = \frac{ns}{n+4s}$. Let $u_{q,v}$ be the solution of

$$L_h^2 u_{q,v} = q|u|^{p-1} v \quad.$$

Then, since $\hat{s} < \frac{n}{4}$, using the first estimate in (3.1.1) and the Hölder inequality, it turns out that

$$\|u_{q,v}\|_s \leq C(n,s) \left\|q|u|^{p-1}\right\|_{\frac{n}{4}} \|v\|_s \quad.$$

Hence the operator $\mathcal{H}_q : v \to u_{q,v}$ acts from $L^s(\mathbb{S}^n)$ to $L^s(\mathbb{S}^n)$ and its norm is less or equal than $C(n,s) \left\|q|u|^{p-1}\right\|_{\frac{n}{4}}$. In particular if $C(n,s) \left\|q|u|^{p-1}\right\|_{\frac{n}{4}} < \frac{1}{2}$, the operator $I - \mathcal{H}_q$ is invertible and the norm of its inverse is less than 2. Now, the equation satisfied by $u$ can be written as

$$(I - \mathcal{H}_q)u = L_h^{-2}(w) \quad.$$

Hence for $s$ fixed and if $\left\|q|u|^{p-1}\right\|_{\frac{n}{4}} < \frac{1}{2}C(n,s)^{-1}$, we deduce

$$\|u\|_s \leq 2C'(n,s) \|w\|_s \leq 2C''(n,s) \|w\|_\infty \quad.$$

Now we use the equation satisfied by $u$ to deduce the estimate of the Proposition. First we take $s > \frac{\frac{n}{4} p\alpha}{\alpha - \frac{n}{4}}$ (this choice implies that $\frac{\alpha s}{s+\alpha p} > \frac{n}{4}$). Then for this $s$ we have (note that $p \leq 2^\# - 1$ allows us to take a $s$ independent of $p$)

$$\|u\|_s \leq C(n,\alpha) \|w\|_\infty \quad .$$

Now

$$\|q|u|^p\|_{\frac{\alpha s}{s+\alpha p}} \leq \|q\|_\alpha \|u\|_s^p \leq \tilde{C}(n,\alpha) \|q\|_\alpha \|w\|_\infty^p \quad ,$$

with $\frac{\alpha s}{s+\alpha p} > \frac{n}{4}$. It follows from the second estimate of (3.1.1) that

$$\|u\|_\infty \leq C(n,\alpha,\|q\|_\alpha) (\|w\|_\infty + \|w\|_\infty^p) \quad .$$

This concludes the proof of the Proposition. □

**Remarks 3.2.** If $\alpha = +\infty$ then $\delta^\infty$ depends only on $n$. As a matter of fact, in this case we choose $s > p\frac{n}{4}$ and then we get

$$\|q|u|^p\|_{\frac{s}{p}} \leq C(n,\|q\|_\infty) \|u\|_s^p \leq C(n,\|q\|_\infty) \|w\|_\infty^p \quad .$$

**Remarks 3.3.** From the arguments in the proof of Proposition 3.1, it follows that we can obtain, for $\left\|q|u|^{p-1}\right\|_{\frac{n}{4}} \leq \delta^\infty$, the estimate

$$\|u\|_\infty \leq C(n,s,\alpha) (\|w\|_s + \|w\|_s^p) \quad ,$$

for all $s > \frac{\frac{n}{4} p\alpha}{\alpha - \frac{n}{4}}$ (and if $\alpha = +\infty$ for all $s > p\frac{n}{4}$ and $C = C(n,s)$).

## §4. A fourth order Sobolev-Aubin inequality.

Following an original idea by Aubin [4], the Sobolev inequality in Theorem 2.1 can be improved when we restrict $u$ to belong to the symmetric class $\mathcal{S}_{2^\#}$ (see notation below). For this new framework, such an improvement is given by Theorem 4.1, to which this section is devoted. In the last steps of the proof (see Lemmas 4.8 and 4.9) are needed some $L^\infty$ estimates proved using Proposition 3.1. The aim of this section is to prove the following Theorem

**Theorem 4.1.** For $q \in (1; 2^\#]$ let $\mathcal{S}_q = \{u \in H_2^2(\mathbb{S}^n) \mid \fint_{\mathbb{S}^n} |u|^q \xi = 0 \text{ for all } \xi \in \mathbb{H}\}$. Then $\exists a_n < 1$ and $\exists q_0 < 2^\#$ so that for all $q_0 \leq q \leq 2^\#$

$$\inf_{\{u \in \mathcal{S}_q \,,\, u \neq 0\}} \frac{a_n \left(\fint_{\mathbb{S}^n} (\Delta_h u)^2 + c_n \fint_{\mathbb{S}^n} |\nabla u|^2\right) + d_n \fint_{\mathbb{S}^n} u^2}{\left(\fint_{\mathbb{S}^n} |u|^q\right)^{\frac{2}{q}}} \geq K_0^{-1} \omega_n^{-\frac{4}{n}} \quad .$$

Before giving the proof of Theorem 4.1, we introduce some further notation.

**Notation 4.2.** For $a < 1$ and $u \in H_2^2(\mathbb{S}^n)$ set

$$E_a[u] = a \left(\fint_{\mathbb{S}^n} (\Delta_h u)^2 + c_n \fint_{\mathbb{S}^n} |\nabla u|^2\right) + d_n \fint_{\mathbb{S}^n} u^2 \quad ;$$

$$J_{a,q}[u] = \frac{E_a[u]}{\left(\fint_{\mathbb{S}^n} |u|^q\right)^{\frac{2}{q}}} \quad \forall u \in H_2^2(\mathbb{S}^n) \setminus \{0\} \quad ;$$

$$\mathcal{M}_{a,q} = \inf_{\{u \in \mathcal{S}_q \,,\, u \neq 0\}} J_{a,q}[u] = \inf_{u \in \mathcal{S}_q^0} E_a[u] \quad ;$$

where
$$\mathcal{S}_q^0 = \left\{ u \in \mathcal{S}_q \text{ such that } \fint_{\mathbb{S}^n} |u|^q = 1 \right\} \quad .$$

We split the proof of Theorem 4.1 in several Lemmas.

**Lemma 4.3.** $\mathcal{M}_{a,q} \leq K_0^{-1} \omega_n^{-\frac{4}{n}}$ for all $0 \leq a < 1$ and for all $2 < q \leq 2^{\#}$. Moreover, for all $2 < q \leq 2^{\#}$, we have

(4.3.1) $$\lim_{a \to 1} \mathcal{M}_{a,q} = K_0^{-1} \omega_n^{-\frac{4}{n}} \quad .$$

**Proof.** We have $\mathcal{M}_{a,q} \leq J_{a,q}[1] = d_n$ and $d_n = K_0^{-1} \omega_n^{-\frac{4}{n}}$ by (1.0.3) and (2.0.1). Hence the first part of the Lemma is proved. Now, according to Hölder's inequality and to Theorem 2.1, for all $2 < q \leq 2^{\#}$

(4.3.2) $$\mathcal{M}_{1,q} \geq K_0^{-1} \omega_n^{-\frac{4}{n}} \quad .$$

Since $a \to \mathcal{M}_{a,q}$ is increasing, $\lim_{a \to 1} \mathcal{M}_{a,q}$ exists. Assume by contradiction that $\lim_{a \to 1} \mathcal{M}_{a,q} < K_0^{-1} \omega_n^{-\frac{4}{n}}$. Setting $\lambda = \lim_{a \to 1} \mathcal{M}_{a,q}$, we have

$$\exists a_0 < 1 \text{ such that } \forall a \in (a_0; 1) \ \exists u_a \in \mathcal{S}_q^0 \quad \text{satisfying} \quad E_a[u_a] \leq \lambda + \frac{l}{4} < K_0^{-1} \omega_n^{-\frac{4}{n}} \quad ,$$

where $l = K_0^{-1} \omega_n^{-\frac{4}{n}} - \lambda$. Hence, $\forall a \in (a_0; 1)$, there exists $u_a \in \mathcal{S}_q^0$ such that

$$E_1[u_a] \leq \lambda + \frac{l}{4} + (1-a) \left( \fint_{\mathbb{S}^n} (\Delta_h u_a)^2 + c_n \fint_{\mathbb{S}^n} |\nabla u_a|^2 \right) \quad .$$

Since the quantities $\fint_{\mathbb{S}^n} (\Delta_h u_a)^2 + c_n \fint_{\mathbb{S}^n} |\nabla u_a|^2$ are uniformly bounded, for $a$ close to 1,

$$E_1[u_a] \leq \lambda + \frac{l}{2} < K_0^{-1} \omega_n^{-\frac{4}{n}} \quad ,$$

contradicting (4.3.2). The second part of the Lemma follows. This concludes the proof of Lemma 4.3. □

**Lemma 4.4.** Let $r \in (2; 2^{\#}]$. For all $\varepsilon > 0$, there exists a $C_\varepsilon$ depending only on $\varepsilon$, so that $\forall u \in \mathcal{S}_r$

$$\left( \fint_{\mathbb{S}^n} |u|^r \right)^{\frac{2}{r}} \leq K_0 \omega_n^{\frac{4}{n}} \left( 2^{-\frac{4}{n}} + \varepsilon \right) \fint_{\mathbb{S}^n} (\Delta_h u)^2 + C_\varepsilon \left( \fint_{\mathbb{S}^n} \left\{ |\nabla u|^2 + u^2 \right\} \right) \quad .$$

**Proof.** This Lemma (analogous to a result of Aubin [4]) is due to Jourdain [29]. □

Let $a_0 \in (2^{-\frac{4}{n}}; 1)$ and fix it in the sequel.

**Lemma 4.5.** For each $a \in [a_0; 1)$ and each $q \in (2; 2^{\#})$, $\inf_{u \in \mathcal{S}_q^0} J_{a,q}[u]$ is achieved by some function $u_{a,q} \in \mathcal{S}_q^0 \cap C^\infty(\mathbb{S}^n)$. Moreover, there exists $\Lambda_{a,q} \in \mathbb{R}^+$ and $\xi_{a,q} \in \mathbb{H}$, $\|\xi_{a,q}\|_\infty = 1$, such that $u_{a,q}$ satisfies

(4.5.1) $$a \left( \Delta_h^2 u_{a,q} + c_n \Delta_h u_{a,q} \right) + d_n u_{a,q} = \mathcal{M}_{a,q} |u_{a,q}|^{q-2} u_{a,q} + \Lambda_{a,q} \xi_{a,q} |u_{a,q}|^{q-2} u_{a,q}$$

and there exists $C > 0$ depending only on $n$ such that for all $a \in [a_0; 1)$ and for all $q \in (2; 2^{\#})$, $\|u_{a,q}\|_2^2 \geq C$.

**Proof.** Using Lemma 4.4 we find the existence of $u_{a,q} \in \mathcal{S}_q^0$, $\Lambda_{a,q} \in \mathbb{R}^+$ and $\xi_{a,q} \in \mathbb{H}$ such that $u_{a,q}$ minimizes $J_{a,q}$ and satisfies the equation

$$a\left(\Delta_h^2 u_{a,q} + c_n \Delta_h u_{a,q}\right) + d_n u_{a,q} = \mathcal{M}_{a,q} |u_{a,q}|^{q-2} u_{a,q} + \Lambda_{a,q} \xi_{a,q} |u_{a,q}|^{q-2} u_{a,q} \quad .$$

Concerning the regularity of $u_{a,q}$, reasoning as in Djadli-Hebey-Ledoux [21] Lemma 2.1, we can prove by a bootstrap argument that for all $s \geq 1$, $u_{a,q} \in L^s(\mathbb{S}^n)$. Then it easily follows that $u_{a,q} \in C^\infty(\mathbb{S}^n)$. Now, since $u_{a,q} \in \mathcal{S}_q^0$, according to Lemma 4.4

$$1 \leq K_0 \omega_n^{\frac{4}{n}} \left(2^{-\frac{4}{n}} + \varepsilon\right) \fint_{\mathbb{S}^n} (\Delta_h u_{a,q})^2 + C_\varepsilon \fint_{\mathbb{S}^n} \left\{|\nabla u_{a,q}|^2 + u_{a,q}^2\right\} \quad .$$

On the other hand by Lemma 4.3, since $u_{a,q}$ is a minimizer

$$a\left(\fint_{\mathbb{S}^n} (\Delta_h u_{a,q})^2 + c_n \fint_{\mathbb{S}^n} |\nabla u_{a,q}|^2\right) + d_n \fint_{\mathbb{S}^n} u_{a,q}^2 \leq K_0^{-1} \omega_n^{-\frac{4}{n}} \quad .$$

Hence

$$C_\varepsilon \|u_{a,q}\|_{H_1^2}^2 \geq 1 - K_0 \omega_n^{\frac{4}{n}} \left(2^{-\frac{4}{n}} + \varepsilon\right) \fint_{\mathbb{S}^n} (\Delta_h u_{a,q})^2$$

$$\geq 1 - K_0 \omega_n^{\frac{4}{n}} \left(\frac{2^{-\frac{4}{n}} + \varepsilon}{a}\right) \left\{K_0^{-1} \omega_n^{-\frac{4}{n}} - c_n a \|u_{a,q}\|_{H_1^2}^2\right\}$$

(note that when $n \geq 5$, $c_n \leq d_n$). Then

$$\left(C_\varepsilon - K_0 c_n \omega_n^{\frac{4}{n}} \left(2^{-\frac{4}{n}} + \varepsilon\right)\right) \|u_{a,q}\|_{H_1^2}^2 \geq 1 - \left(\frac{2^{-\frac{4}{n}} + \varepsilon}{a_0}\right) > 0$$

for $\varepsilon$ small. Now, since $H_2^2 \hookrightarrow H_1^2 \hookrightarrow L^2$, the first embedding being compact and the second one being continuous, and using an interpolation type inequality of Lions [36], we get: for all $\delta > 0$ $\exists C_\delta > 0$ such that for all $u \in H_2^2(\mathbb{S}^n)$

$$\|u_{a,q}\|_{H_1^2}^2 \leq \delta \|u_{a,q}\|_{H_2^2}^2 + C_\delta \|u_{a,q}\|_2^2 \quad .$$

It follows that for $\delta$ small enough

$$C_\delta \|u_{a,q}\|_2^2 \geq C(n) - \delta \|u_{a,q}\|_{H_2^2}^2 \geq C(n) - \frac{\delta}{a} K_0^{-1} \omega_n^{-\frac{4}{n}} \geq C(n) - \delta 2^{\frac{4}{n}} K_0^{-1} \omega_n^{-\frac{4}{n}} \geq C'(n) > 0 \quad ,$$

and we get the result. This concludes the proof of Lemma 4.5. $\square$

**Lemma 4.6.** *For all $a \in [a_0; 1)$ and all $q \in (2; 2^\#)$, $\Lambda_{a,q} \leq C(n)$, where $C(n)$ is a constant depending only on $n$.*

**Proof.** Multiplying (4.5.1) by $u_{a,q} \xi_{a,q}$ and integrating over $\mathbb{S}^n$, we get

$$a\left(\fint_{\mathbb{S}^n} \{\Delta_h^2 u_{a,q} + c_n \Delta_h u_{a,q}\} u_{a,q} \xi_{a,q}\right) + d_n \fint_{\mathbb{S}^n} u_{a,q}^2 \xi_{a,q} = \Lambda_{a,q} \fint_{\mathbb{S}^n} |u_{a,q}|^q \xi_{a,q}^2 \quad .$$

Integrating by parts and taking into account $\Delta_h \xi_{a,q} = n \xi_{a,q}$, we have

$$\fint_{\mathbb{S}^n} (\Delta_h^2 u_{a,q}) u_{a,q} \xi_{a,q} = \fint_{\mathbb{S}^n} (\Delta_h u_{a,q})^2 \xi_{a,q} + n \fint_{\mathbb{S}^n} (\Delta_h u_{a,q}) u_{a,q} \xi_{a,q} - 2 \fint_{\mathbb{S}^n} (\Delta_h u_{a,q}) \nabla u_{a,q} \nabla \xi_{a,q} \quad ;$$

moreover

$$\fint_{\mathbb{S}^n} \Delta_h u_{a,q} \nabla u_{a,q} \nabla \xi_{a,q} = \frac{n}{2} \fint_{\mathbb{S}^n} |\nabla u_{a,q}|^2 \xi_{a,q} \quad .$$

It follows that

$$\Lambda_{a,q}\fint_{\mathbb{S}^n} |u_{a,q}|^q \xi_{a,q}^2 = a\left(\fint_{\mathbb{S}^n} (\Delta_h u_{a,q})^2 \xi_{a,q} + (c_n+n)\fint_{\mathbb{S}^n} (\Delta_h u_{a,q})\, u_{a,q}\xi_{a,q}\right)$$
$$- na\fint_{\mathbb{S}^n} |\nabla u_{a,q}|^2 \xi_{a,q} + d_n \fint_{\mathbb{S}^n} u_{a,q}^2 \xi_{a,q} \quad.$$

Independently

$$\fint_{\mathbb{S}^n} (\Delta_h u_{a,q})\, u_{a,q}\xi_{a,q} = \fint_{\mathbb{S}^n} |\nabla u_{a,q}|^2 \xi_{a,q} + \frac{n}{2}\fint_{\mathbb{S}^n} (u_{a,q})^2 \xi_{a,q}$$

and hence,

$$\Lambda_{a,q}\fint_{\mathbb{S}^n} |u_{a,q}|^q \xi_{a,q}^2 = a\left(\fint_{\mathbb{S}^n} (\Delta_h u_{a,q})^2 \xi_{a,q} + c_n\fint_{\mathbb{S}^n} \left(|\nabla u_{a,q}|^2 \xi_{a,q}\right)\right)$$
$$+ \left(\frac{n^2 a}{2} + c_n a \frac{n}{2} + d_n\right)\fint_{\mathbb{S}^n} u_{a,q}^2 \xi_{a,q} \quad.$$

It follows that

$$\Lambda_{a,q}\fint_{\mathbb{S}^n} |u_{a,q}|^q \xi_{a,q}^2 \leq C(n) E_a[u_{a,q}] \leq C(n) \quad.$$

Using Lemma 4.5 we have $\|u_{a,q}\|_2^2 \geq C$, hence by the Hölder inequality there exists a constant $C' > 0$ such that for all $a$ and $q$, $\fint_{\mathbb{S}^n} |u_{a,q}|^q \xi_{a,q}^2 \geq C'$. It follows that $\Lambda_{a,q} \leq C(n)$. This concludes the proof of Lemma 4.6. □

We now consider increasing sequences $(a_k)$ and $(q_k)$ such that

$$a_k < 1 \quad,\quad a_k \to 1 \quad;\quad 2 < q_k < 2^\# \quad,\quad q_k \to 2^\# \quad.$$

By Lemma 4.5 $u_k$, the solution of (4.5.1) associated to $a_k$ and $q_k$, satisfies

$$a_k\left(\Delta_h^2 u_k + c_n \Delta_h u_k\right) + d_n u_k = \mathcal{M}_k |u_k|^{q_k-2} u_k + \Lambda_k \xi_k |u_k|^{q_k-2} u_k \quad.$$

Passing to a subsequence of $u_k$ if necessary, and still denoting the subsequence by $u_k$, we will establish the following properties of $(u_k)$ :

$$\begin{cases} (a) & u_k \to 1 \text{ weakly in } H_2^2(\mathbb{S}^n) \text{ and } u_k \to 1 \text{ in } L^{2^\#}(\mathbb{S}^n) \quad; \\ (b) & \text{there exists } C(n) \text{ such that } \|u_k\|_\infty \leq C(n) \quad; \\ (c) & \|u_k - 1\|_\infty = o(1) \text{ as } k \to +\infty \quad. \end{cases}$$

One has first to remark that $\lim_{k\to+\infty} \mathcal{M}_k = K_0^{-1}\omega_n^{-\frac{4}{n}}$. As a matter of fact, since $(\mathcal{M}_k)$ is a bounded sequence (see Lemma 4.3), we can extract a subsequence such that $\lim_{k\to+\infty} \mathcal{M}_k$ exists. Assuming that $\lim_{k\to+\infty} \mathcal{M}_k = l < K_0^{-1}\omega_n^{-\frac{4}{n}}$ leads to a contradiction, as in the proof of Lemma 4.3.

Now, since $E_{a_k}[u_k]$ is bounded, $(u_k)$ is bounded in $H_2^2(\mathbb{S}^n)$; hence, some subsequence of $(u_k)$ converges weakly in $H_2^2(\mathbb{S}^n)$ and almost everywhere on $\mathbb{S}^n$ to a function $u \in C^\infty(\mathbb{S}^n)$ which satisfies

(4.6.1) $$\Delta_h^2 u + c_n \Delta_h u + d_n u = K_0^{-1}\omega_n^{-\frac{4}{n}} |u|^{\frac{8}{n-4}} u + \Lambda \xi_\infty |u|^{\frac{8}{n-4}} u \quad.$$

Here $\Lambda \geq 0$ is the limit of $(\Lambda_k)$, and $\xi_\infty \in \mathbb{H}$, $\|\xi_\infty\|_\infty = 1$, is the limit of $(\xi_k)$ (recall that $\mathbb{H}$ is finite dimensional).

**Lemma 4.7.** $u_k \to 1$ weakly in $H_2^2(\mathbb{S}^n)$ and strongly in $L^{2^\#}(\mathbb{S}^n)$.

**Proof.** First we remark that according to Lemma 4.5, $u$ is not identically 0. Using Theorem 1.1, we can write

$$\Lambda \fint_{\mathbb{S}^n} |\nabla \xi_\infty|^2 |u|^{2^\#} = 0 \quad,$$

so it follows that $\Lambda = 0$. We multiply (4.6.1) by $u$ and we integrate to obtain

$$\fint_{\mathbb{S}^n} (\Delta u)^2 + c_n \fint_{\mathbb{S}^n} |\nabla u|^2 + d_n \fint_{\mathbb{S}^n} u^2 = K_0^{-1} \omega_n^{-\frac{4}{n}} \fint_{\mathbb{S}^n} |u|^{2^\#} \quad.$$

According to Theorem 2.1, this gives

$$K_0^{-1} \omega_n^{-\frac{4}{n}} \fint_{\mathbb{S}^n} |u|^{2^\#} \geq K_0^{-1} \omega_n^{-\frac{4}{n}} \left( \fint_{\mathbb{S}^n} |u|^{2^\#} \right)^{\frac{2}{2^\#}} \quad,$$

so it follows that $\fint_{\mathbb{S}^n} |u|^{2^\#} \geq 1$. Independently, since $(u_k)$ is bounded in $H_2^2(\mathbb{S}^n)$

$$K_0^{-1} \omega_n^{-\frac{4}{n}} \left( \fint_{\mathbb{S}^n} |u_k|^{2^\#} \right)^{\frac{2}{2^\#}} \leq E_1[u_k] \leq \mathcal{M}_k + (1 - a_k) \left( \fint_{\mathbb{S}^n} (\Delta u_k)^2 + c_n \fint_{\mathbb{S}^n} |\nabla u_k|^2 \right) \to K_0^{-1} \omega_n^{-\frac{4}{n}} \quad.$$

Then,

$$1 \leq \fint_{\mathbb{S}^n} |u|^{2^\#} \leq \liminf_{k \to \infty} \fint_{\mathbb{S}^n} |u_k|^{2^\#} \leq \limsup_{k \to \infty} \fint_{\mathbb{S}^n} |u_k|^{2^\#} \leq 1 \quad.$$

Hence, $\fint_{\mathbb{S}^n} |u|^{2^\#} = 1$ and $\lim_{k \to +\infty} \fint_{\mathbb{S}^n} |u_k|^{2^\#} = 1$. Then, applying a result of Brézis-Lieb [12] (see also Kavian [30]), we have $\lim_{k \to +\infty} \|u - u_k\|_{2^\#}^{2^\#} = \lim_{k \to +\infty} \|u_k\|_{2^\#}^{2^\#} - \|u\|_{2^\#}^{2^\#} = 0$. This means that $u_k \to u$ in $L^{2^\#}(\mathbb{S}^n)$. It remains to prove that $u \equiv 1$. First we show that $\fint_{\mathbb{S}^n} |u|^{2^\#} \xi = 0$ for all $\xi \in \mathbb{H}$. We have $\fint_{\mathbb{S}^n} |u_k|^{q_k - 1} |u| \xi \to \fint_{\mathbb{S}^n} |u|^{2^\#} \xi$ as $k \to +\infty$ since $|u_k|^{q_k - 1} \to |u|^{2^\# - 1}$ weakly in $L^{\frac{2^\#}{2^\# - 1}}(\mathbb{S}^n)$. On the other hand, since $u_k \in \mathcal{S}_{q_k}$

$$\left| \fint_{\mathbb{S}^n} |u_k|^{q_k - 1} |u| \xi \right| = \left| \fint_{\mathbb{S}^n} |u_k|^{q_k - 1} (|u| - |u_k|) \xi \right| \leq \|u_k\|_{2^\#}^{q_k - 1} \|u - u_k\|_r \leq C \|u - u_k\|_r \quad,$$

where $r = \frac{2^\#}{2^\# - q_k + 1} \leq 2^\#$. Thus

$$\fint_{\mathbb{S}^n} |u|^{2^\#} \xi = \lim_{k \to +\infty} \fint_{\mathbb{S}^n} |u_k|^{q_k - 1} |u| \xi = 0 \quad.$$

This proves that $u \in \mathcal{S}_{2^\#}^0$. Hence $u$ is a solution of (4.6.1) with $u \in \mathcal{S}_{2^\#}^0$. Since such a solution is unique, see Edmunds-Fortunato-Janelli [25] and Lions [37], we have $u \equiv 1$ and this concludes the proof. $\square$

**Lemma 4.8.** $(u_k)$ is bounded in $L^\infty(\mathbb{S}^n)$.

**Proof.** For $k_0 \in \mathbb{N}$, let

$$\Omega_{k_0} = \{ x \in \mathbb{S}^n \mid |u_k(x)| < 2 \quad \text{for all } k \geq k_0 \} \quad.$$

Since $u_k \to 1$ in $L^{2^\#}(\mathbb{S}^n)$ it is clear that $|\mathbb{S}^n \setminus \Omega_{k_0}| = measure(\mathbb{S}^n \setminus \Omega_{k_0}) \to 0$ as $k_0 \to +\infty$. We can write

$$\left( \Delta_h + \frac{c_n}{2} \right)^2 u_k = b_k |u_k|^{q_k - 2} + w_k \quad,$$

where

$$b_k = \left(\left(\frac{\mathcal{M}_k + \Lambda_k \xi_k}{a_k}\right) + \left(\frac{c_n^2}{4} - \frac{d_n}{a_k}\right)\frac{1}{|u_k|^{q_k-2}}\right)1_{\mathbb{S}^n\setminus\Omega_{k_0}} \quad,$$

$$w_k = \left(\left(\frac{\mathcal{M}_k + \Lambda_k \xi_k}{a_k}\right)|u_k|^{q_k-2} + \left(\frac{c_n^2}{4} - \frac{d_n}{a_k}\right)\right)u_k.1_{\Omega_{k_0}} \quad.$$

According to lemmas 4.3 and 4.6, there exists a constant $C > 0$ depending only on $n$ such that $\forall k \in \mathbb{N}$

(4.8.1) $$\|b_k\|_\infty \leq C \quad \text{and} \quad \|w_k\|_\infty \leq C \quad,$$

and

(4.8.2) $$\left\|b_k |u_k|^{q_k-2}\right\|_{\frac{n}{4}} \leq C\left(\|u_k - 1\|_{2^\#} + |\mathbb{S}^n \setminus \Omega_{k_0}|^{\frac{4}{n}} + |\mathbb{S}^n \setminus \Omega_{k_0}|^{\frac{n-4}{2n}}\right) \quad.$$

Using (4.8.1), the fact that $u_k \to 1$ in $L^{2^\#}(\mathbb{S}^n)$ as $k \to +\infty$ and the fact that $|\mathbb{S}^n \setminus \Omega_{k_0}| \to 0$ as $k_0 \to +\infty$, it follows that

$$\left\|b_k |u_k|^{q_k-2}\right\|_{\frac{n}{4}} \to 0 \text{ as } k \to +\infty \quad.$$

Then, according to Proposition 3.1, Remark 3.2 and (4.8.1), there exists a constant $C$ depending only on $n$ such that for all $k$ large

$$\|u_k\|_\infty \leq C \quad.$$

This concludes the proof of the Lemma. □

**Lemma 4.9.** $\|u_k - 1\|_\infty = o(1)$ as $k \to +\infty$.

**Proof.** Fix $\eta > 0$ and for $k_0 \in \mathbb{N}$, let

$$\Omega_{k_0} = \{x \in \mathbb{S}^n \mid |u_k(x) - 1| < \eta \quad \text{for all } k \geq k_0\} \quad.$$

As in the proof of the previous Lemma, $|\mathbb{S}^n \setminus \Omega_{k_0}| \to 0$ as $k_0 \to +\infty$. We have

$$\left(\Delta_h + \frac{c_n}{2}\right)^2 (u_k - 1) = b_k(u_k - 1) + w_k \quad,$$

where

$$b_k = \left(\left(\frac{\mathcal{M}_k + \Lambda_k \xi_k}{a_k}\right)|u_k|^{q_k-2} + \left(\frac{c_n^2}{4} - \frac{d_n}{a_k}\right)\right)1_{\mathbb{S}^n\setminus\Omega_{k_0}} \quad,$$

$$w_k = \left(\left(\frac{\mathcal{M}_k + \Lambda_k \xi_k}{a_k}\right)|u_k|^{q_k-2} + \left(\frac{c_n^2}{4} - \frac{d_n}{a_k}\right)\right)(u_k - 1)1_{\Omega_{k_0}} + \left(\left(\frac{\mathcal{M}_k + \Lambda_k \xi_k}{a_k}\right)|u_k|^{q_k-2} - \frac{d_n}{a_k}\right) \quad.$$

As before there exists $C > 0$ depending only on $n$ such that for all $k \in \mathbb{N}$

$$\|b_k\|_\infty \leq C \quad \text{and} \quad \|w_k\|_\infty \leq C \quad,$$
$$\|b_k\|_{\frac{n}{4}} \leq C |\mathbb{S}^n \setminus \Omega_{k_0}|^{\frac{4}{n}} \quad.$$

Consider $\delta^\infty$ given by Remark 3.2, and $s > \frac{n}{4}$; since $|\mathbb{S}^n \setminus \Omega_{k_0}| \to 0$ as $k_0 \to +\infty$, using Remark 3.3, there exists $C$, depending only on $n$, and there exists $k_1$ such that for all $k \geq k_1$

$$\|u_k - 1\|_\infty \leq C \|w_k\|_s \quad.$$

Clearly, using Lemma 4.8 and the fact that $u_k \to 1$ almost everywhere, there exists $k_2 \in \mathbb{N}$ such that for all $k \geq k_2$

$$\|w_k\|_s \leq C\eta \quad,$$

where $C$ is independent of $k$ and $\eta$. Then we have the Lemma. □

We can now prove Theorem 4.1

**Proof of Theorem 4.1 :** We write $u_k = 1 + \alpha_k + h_k + \Psi_k$ where $\alpha_k \in \mathbb{R}$, $h_k \in \mathbb{H}$ and where $\Psi_k$ is orthogonal to $\mathbb{R} \oplus \mathbb{H}$. Clearly,

$$E_{a_k}[u_k] \geq d_n + a_k \left( \fint_{\mathbb{S}^n} (\Delta_h \Psi_k)^2 + c_n \fint_{\mathbb{S}^n} |\nabla \Psi_k|^2 \right) + 2\alpha_k d_n + d_n \fint_{\mathbb{S}^n} \Psi_k^2 \ .$$

We want to estimate $\alpha_k$. For this, we know that (using Lemma 4.9)

$$0 = \fint_{\mathbb{S}^n} (1 + \alpha_k + h_k + \Psi_k)^{q_k} h_k = q_k \fint_{\mathbb{S}^n} h_k^2 + o(1)\alpha_k \left( \fint_{\mathbb{S}^n} h_k^2 \right)^{\frac{1}{2}} + o(1)\fint_{\mathbb{S}^n} h_k^2 + o(1) \left( \fint_{\mathbb{S}^n} h_k^2 \right)^{\frac{1}{2}} \left( \fint_{\mathbb{S}^n} \Psi_k^2 \right)^{\frac{1}{2}}$$

and this implies, since $\alpha_k \leq 0$ (recall that $u_k \in \mathcal{S}^0_{q_k}$)

$$\fint_{\mathbb{S}^n} h_k^2 = o(1) \left( \alpha_k^2 + \fint_{\mathbb{S}^n} \Psi_k^2 \right) \ .$$

Secondly, we have

$$1 = \fint_{\mathbb{S}^n} (1 + \alpha_k + h_k + \Psi_k)^{q_k} = 1 + q_k \alpha_k + \frac{q_k(q_k-1)}{2}\alpha_k^2 + \frac{q_k(q_k-1)}{2}\fint_{\mathbb{S}^n} \Psi_k^2 + o(1)\alpha_k^2 + o(1)\fint_{\mathbb{S}^n} \Psi_k^2 \ ,$$

and this gives (using Lemma 4.7 to see that $\alpha_k \to 0$)

$$\alpha_k = -\frac{q_k-1}{2}\fint_{\mathbb{S}^n} \Psi_k^2 + o(1)\fint_{\mathbb{S}^n} \Psi_k^2 \ .$$

Now, using the Courant-Fischer characterization of the eigenvalues and Theorem 2.2

$$E_{a_k}[u_k] \geq d_n + a_k \left( \fint_{\mathbb{S}^n} (\Delta_h \Psi_k)^2 + c_n \fint_{\mathbb{S}^n} |\nabla \Psi_k|^2 \right) - (d_n(q_k-2) + o(1)) \fint_{\mathbb{S}^n} \Psi_k^2$$

$$\geq d_n + \left( 4(n+1)^2 + 2c_n(n+1) - \frac{8d_n}{n-4} + (a_k - 1)(n+1)(4(n+1) + 2c_n) + o(1) \right) \fint_{\mathbb{S}^n} \Psi_k^2$$

$$\geq d_n + \frac{1}{2}\left( n^3 + 6n^2 + 8n + \beta_k \right) \fint_{\mathbb{S}^n} \Psi_k^2$$

where $\beta_k \to 0$ when $k$ goes to $+\infty$. This implies that for $k$ large

$$E_{a_k}[u_k] \geq d_n$$

and this proves that $\exists\, a_n < 1$ and $\exists\, q_0 < 2^\#$ so that

$$\inf_{\{u \in \mathcal{S}_q\ ,\ u \neq 0\}} \frac{a_n \left( \fint_{\mathbb{S}^n} (\Delta_h u)^2 + c_n \fint_{\mathbb{S}^n} |\nabla u|^2 \right) + d_n \fint_{\mathbb{S}^n} u^2}{\left( \fint_{\mathbb{S}^n} |u|^q \right)^{\frac{2}{q}}} \geq K_0^{-1} \omega_n^{-\frac{4}{n}}$$

for all $q_0 \leq q < 2^\#$. It remains to prove that this is also true for $q = 2^\#$. Denote

$$\mathcal{M}_{2^\#} = \inf_{\{u \in \mathcal{S}_{2^\#}\ ,\ u \neq 0\}} \frac{a_n \left( \fint_{\mathbb{S}^n} (\Delta_h u)^2 + c_n \fint_{\mathbb{S}^n} |\nabla u|^2 \right) + d_n \fint_{\mathbb{S}^n} u^2}{\left( \fint_{\mathbb{S}^n} |u|^{2^\#} \right)^{\frac{2}{2^\#}}}$$

and consider $\varepsilon > 0$; then there exists $u \in \mathcal{S}_{2^\#}$ such that $J_{a_n,2^\#}[u] \leq \mathcal{M}_{2^\#} + \varepsilon$. Since

$$\lim_{q \to 2^\#} \left( \fint_{\mathbb{S}^n} |u|^q \right)^{\frac{2}{q}} = \left( \fint_{\mathbb{S}^n} |u|^{2^\#} \right)^{\frac{2}{2^\#}} \quad,$$

it is clear that

$$K_0^{-1} \omega_n^{-\frac{4}{n}} \leq \limsup_{q \to 2^\#} J_{a_n,q}[u] = J_{a_n,2^\#}[u] \leq \mathcal{M}_{2^\#} + \varepsilon \quad,$$

and then $K_0^{-1} \omega_n^{-\frac{4}{n}} \leq \mathcal{M}_{2^\#}$. This concludes the proof of the Theorem. $\square$

## §5. Basic estimates.

Given $P \in \mathbb{S}^n$ and $t \in [1; +\infty)$ let $\varphi_{P,t}$ be the conformal transformation on $\mathbb{S}^n$ defined in section 1. For brevity we will denote $\varphi_{P,t}$ as $\varphi_p$ and $f_{P,t} = f \circ \varphi_{P,t}$ as $f_p$. We set $T_\varphi u = (u \circ \varphi) |det\, d\varphi|^{\frac{n-4}{2n}}$. Clearly, $\varphi_{P,t}$ acts on the set of conformal metrics $g = u^{\frac{4}{n-4}} h$ by

$$\varphi_{P,t}^\star g = \left( T_{\varphi_{P,t}} u \right)^{\frac{4}{n-4}} h \quad.$$

Consequently, equation $(E)$ is transformed by $\varphi_{P,t}$ into

$$\Delta_h^2 \left( T_{\varphi_{P,t}} u \right) + c_n \Delta_h \left( T_{\varphi_{P,t}} u \right) + d_n \left( T_{\varphi_{P,t}} u \right) = \frac{n-4}{2} \left( f \circ \varphi_{P,t} \right) \left( T_{\varphi_{P,t}} u \right)^{\frac{n+4}{n-4}} \quad.$$

Roughly, by means of the improved Sobolev inequality in Theorem 4.1, for $f$ close to a constant, it is possible to find minima of $\dfrac{E_1[u]}{\left( \fint_{\mathbb{S}^n} f|u|^{2^\#} \right)^{\frac{2}{2^\#}}}$ constrained on $\mathcal{S}_{2^\#}$, leading to solutions of

$$(E_p) \qquad P_h^n u = \left( \frac{n-4}{2} f_p - \Lambda_p \xi_p \right) |u|^{2^\# - 2} u \quad,$$

where $\Lambda_p$ is a Lagrange multiplier. Using this finite dimensional reduction, solving problem $(E)$ amounts to finding $p_0$ for which $\Lambda_{p_0} = 0$. The same kind of strategy has been used by Chang-Yang [19].

Let $t > 1$ and $P \in \mathbb{S}^n$. Considering $f_p$ and $2 < q \leq 2^\#$, we define

$$\mathcal{M}_q = \inf_{\left\{ u \in \mathcal{S}_q | \fint_{\mathbb{S}^n} f_p |u|^q = 1 \right\}} E_1[u] \quad;$$

it is now classical that for $q < 2^\#$ there exist $u_q \in C^\infty(\mathbb{S}^n) \cap \mathcal{S}_q$, $\Lambda_q \in \mathbb{R}^+$ and $\xi_q \in \mathbb{H}$, $\|\xi_q\|_\infty = 1$, which satisfy the equation

$$\Delta_h^2 u_q + c_n \Delta_h u_q + d_n u_q = \left( \mathcal{M}_q f_p - \Lambda_q \xi_q \right) |u_q|^{q-2} u_q \quad.$$

Clearly $(u_q)$ is bounded in $H_2^2(\mathbb{S}^n)$ and it follows that $(u_q)$ converges weakly in $H_2^2(\mathbb{S}^n)$ to some $\hat{u}_p \in H_2^2(\mathbb{S}^n)$ satisfying

$$\Delta_h^2 \hat{u}_p + c_n \Delta_h \hat{u}_p + d_n \hat{u}_p = \left( \mathcal{M}_p f_p - \hat{\Lambda}_p \hat{\xi}_p \right) |\hat{u}_p|^{2^\# - 2} \hat{u}_p \quad,$$

where $\hat{\Lambda}_p \in \mathbb{R}^+$ and $\hat{\xi}_p \in \mathbb{H}$, $\left\| \hat{\xi}_p \right\|_\infty = 1$. We claim that $\hat{u}_p$ is not identically 0, and that $\hat{u}_p \in \mathcal{S}_{2^\#}$. The proof of this fact relies on the following claim : if $\varepsilon_f = \left\| f - \frac{n(n^2-4)}{8} \right\|_\infty$ is small enough (this smallness being independent of $t$ and $P$), then $(u_q)_q$ is bounded in $L^\infty(\mathbb{S}^n)$. To show this, one can prove that

$$\begin{cases} 1) \fint_{\mathbb{S}^n} (\Delta_h u_q)^2 + c_n \fint_{\mathbb{S}^n} |\nabla u_q|^2 = O\left(\|f_p - f(P)\|_\infty\right) \quad, \\ \fint_{\mathbb{S}^n} u_q^2 - \left(\frac{n(n^2-4)}{8f(P)}\right)^{\frac{n-4}{4}} = O\left(\|f_p - f(P)\|_\infty\right) \quad; \\ 2) \left\| u_q - \left(\frac{n(n^2-4)}{8f(P)}\right)^{\frac{n-4}{8}} \right\|_{H_2^2}^2 = O\left(\|f_p - f(P)\|_\infty\right) \quad; \\ 3) (\Lambda_q) \text{ is bounded} \quad. \end{cases}$$

reasoning as in Lemmas 5.2, 5.3 and 5.4.

Up to a subsequence, $(u_q)_q$ converges to some $\hat{u}_p$ almost everywhere and $(u_q)_q$ is bounded in $L^\infty(\mathbb{S}^n)$. It follows easily that $(u_q)$ converges to $\hat{u}_p$ in $L^{2^\#}(\mathbb{S}^n)$, where $\hat{u}_p \not\equiv 0$, $\hat{u}_p \in \mathcal{S}_{2^\#}$ and $\hat{u}_p$ realizes $\mathcal{M}_{2^\#}$.

Regarding equation $(E_p)$, we provide uniform estimates on the solutions $u_p$, uniform with respect to $p$. These estimates will be needed later to study the reduced finite dimensional problem.

**Lemma 5.1.** *Denote*

$$\mathcal{M}_p = \inf_{\left\{u \in \mathcal{S}_{2^\#} \;\middle|\; \fint_{\mathbb{S}^n} f_p|u|^{2^\#} = 1\right\}} E_1[u] \quad.$$

*Then for $\varepsilon > 0$ small enough, if $f$ satisfies*

$$\varepsilon_f = \left\| f - \frac{n(n^2-4)}{8} \right\|_\infty \leq \varepsilon \quad,$$

*there exist $u_p \in C^\infty(\mathbb{S}^n) \cap \mathcal{S}_{2^\#}$ with $\fint_{\mathbb{S}^n} u_p \geq 0$, $\Lambda_p \in \mathbb{R}^+$ and $\xi_p \in \mathbb{H}$, $\|\xi_p\|_\infty = 1$, which satisfy the equation*

$$(5.1.1) \qquad \Delta_h^2 u_p + c_n \Delta_h u_p + d_n u_p = \left(\frac{n-4}{2} f_p - \Lambda_p \xi_p\right) |u_p|^{2^\# - 2} u_p \quad.$$

*Furthermore, $(u_p)_p$ is bounded in $H_2^2(\mathbb{S}^n)$ and we have*

$$(5.1.2) \qquad \frac{n-4}{2} \left( \fint_{\mathbb{S}^n} f_p |u_p|^{2^\#} \right)^{\frac{4}{n}} = \mathcal{M}_p \leq d_n \left( \fint_{\mathbb{S}^n} f_p \right)^{-\frac{2}{2^\#}} \quad.$$

**Proof.** As we have just seen, there exist $\hat{u}_p \in C^\infty(\mathbb{S}^n) \cap \mathcal{S}_{2^\#}$, $\hat{\Lambda}_p \in \mathbb{R}^+$ and $\hat{\xi}_p \in \mathbb{H}$, $\|\hat{\xi}_p\|_\infty = 1$, which satisfy the equation

$$\Delta_h^2 \hat{u}_p + c_n \Delta_h \hat{u}_p + d_n \hat{u}_p = \left(\mathcal{M}_p f_p - \hat{\Lambda}_p \hat{\xi}_p\right) |\hat{u}_p|^{2^\# - 2} \hat{u}_p \quad.$$

Hence, by a simple renormalization there exist $u_p \in C^\infty(\mathbb{S}^n) \cap \mathcal{S}_{2^\#}$ with $\fint_{\mathbb{S}^n} u_p \geq 0$, $\Lambda_p \in \mathbb{R}^+$ and $\xi_p \in \mathbb{H}$, $\|\xi_p\|_\infty = 1$, which satisfy equation (5.1.1). Clearly, with this renormalization it is $\frac{n-4}{2} \left( \fint_{\mathbb{S}^n} f_p |u_p|^{2^\#} \right)^{\frac{4}{n}} = \mathcal{M}_p$. Moreover, we have

$$(5.1.3) \qquad \mathcal{M}_p \leq E_1\left[\left(\fint_{\mathbb{S}^n} f_p\right)^{-\frac{1}{2^\#}}\right] = d_n \left(\fint_{\mathbb{S}^n} f_p\right)^{-\frac{2}{2^\#}} \quad.$$

Now, using the estimate $\left(\fint_{\mathbb{S}^n} f_p |u_p|^{2^\#}\right)^{\frac{4}{n}} \leq \frac{2}{n-4} d_n \left(\fint_{\mathbb{S}^n} f_p\right)^{-\frac{2}{2^\#}}$, it is clear that for $\varepsilon$ small enough, $(u_p)$ is bounded in $L^{2^\#}(\mathbb{S}^n)$. Since

$$\fint_{\mathbb{S}^n} (\Delta_h u_p)^2 + c_n \fint_{\mathbb{S}^n} |\nabla u_p|^2 + d_n \fint_{\mathbb{S}^n} u_p^2 = \frac{n-4}{2} \fint_{\mathbb{S}^n} f_p |u_p|^{2^\#} \quad,$$

$(u_p)_p$ is also bounded in $H_2^2(\mathbb{S}^n)$. This concludes the proof of Lemma 5.1. $\square$

**Lemma 5.2.** *The following estimates hold, uniformly in $p$*

$$\begin{cases} \fint_{\mathbb{S}^n} (\Delta_h u_p)^2 + c_n \fint_{\mathbb{S}^n} |\nabla u_p|^2 = O\left(\|f_p - f(P)\|_\infty\right) & ; \\ \fint_{\mathbb{S}^n} u_p^2 - \left(\frac{n(n^2-4)}{8f(P)}\right)^{\frac{n-4}{4}} = O\left(\|f_p - f(P)\|_\infty\right) & . \end{cases}$$

**Proof.** We have by (5.1.2)

$$E_1[u_p] = \frac{n-4}{2} \fint_{\mathbb{S}^n} f_p |u_p|^{2^\#} = \left(\frac{2}{n-4}\right)^{\frac{n}{4}-1} \mathcal{M}_p^{\frac{n}{22^\#}} \mathcal{M}_p \leq \left(\frac{2}{n-4}\right)^{\frac{n}{4}-1} \mathcal{M}_p^{\frac{n}{22^\#}} d_n \left(\fint_{\mathbb{S}^n} f_p\right)^{-\frac{2}{2^\#}}.$$

By Theorem 4.1 applied with $q = 2^\#$, $d_n \left(\fint_{\mathbb{S}^n} |u_p|^{2^\#}\right)^{\frac{2}{2^\#}} \leq a_n \left(\fint_{\mathbb{S}^n} (\Delta_h u_p)^2 + c_n \fint_{\mathbb{S}^n} |\nabla u_p|^2\right) + d_n \fint_{\mathbb{S}^n} u_p^2$.
It follows that

$$(1 - a_n) \left(\fint_{\mathbb{S}^n} (\Delta_h u_p)^2 + c_n \fint_{\mathbb{S}^n} |\nabla u_p|^2\right)$$

$$\leq \left(\frac{2}{n-4}\right)^{\frac{n}{4}-1} \mathcal{M}_p^{\frac{n}{22^\#}} d_n \left(\fint_{\mathbb{S}^n} f_p\right)^{-\frac{2}{2^\#}} - d_n \left(f(P)^{-1} \fint_{\mathbb{S}^n} f(P) |u_p|^{2^\#}\right)^{\frac{2}{2^\#}}$$

$$\leq \left(\frac{2}{n-4}\right)^{\frac{n}{4}-1} \mathcal{M}_p^{\frac{n}{22^\#}} \left(d_n \left(\fint_{\mathbb{S}^n} f_p\right)^{-\frac{2}{2^\#}} - d_n (f(P))^{-\frac{2}{2^\#}}\right) + \left(\frac{2}{n-4}\right)^{\frac{n}{4}-1} \mathcal{M}_p^{\frac{n}{22^\#}} d_n (f(P))^{-\frac{2}{2^\#}}$$

$$- d_n \left(f(P)^{-1} \fint_{\mathbb{S}^n} f(P) |u_p|^{2^\#}\right)^{\frac{2}{2^\#}}$$

$$\leq \left(\frac{2}{n-4}\right)^{\frac{n}{4}-1} \mathcal{M}_p^{\frac{n}{22^\#}} d_n \left(\left(\fint_{\mathbb{S}^n} f_p\right)^{-\frac{2}{2^\#}} - (f(P))^{-\frac{2}{2^\#}}\right)$$

$$+ \frac{d_n}{(f(P))^{\frac{2}{2^\#}}} \left(\left(\fint_{\mathbb{S}^n} f_p |u_p|^{2^\#}\right)^{\frac{2}{2^\#}} - \left(\fint_{\mathbb{S}^n} f(P) |u_p|^{2^\#}\right)^{\frac{2}{2^\#}}\right) \quad.$$

With same computations using Lemma 5.1, we get $\fint_{\mathbb{S}^n} (\Delta_h u_p)^2 + c_n \fint_{\mathbb{S}^n} |\nabla u_p|^2 = O\left(\|f_p - f(P)\|_\infty\right)$. Now

$$E_1[u_p] - d_n \left(\frac{n(n^2-4)}{8f(P)}\right)^{\frac{n-4}{4}} = \left(\frac{2}{n-4}\right)^{\frac{n}{4}-1} \mathcal{M}_p^{\frac{n}{4}} - d_n \left(\frac{n(n^2-4)}{8f(P)}\right)^{\frac{n-4}{4}}$$

$$= \left(\frac{2}{n-4}\right)^{\frac{n}{4}-1} \left(\mathcal{M}_p^{\frac{n}{4}} - d_n \left(\frac{d_n}{f(P)}\right)^{\frac{n-4}{4}}\right) \leq \left(\frac{2}{n-4}\right)^{\frac{n}{4}-1} d_n^{\frac{n}{4}} \left(\left(\fint_{\mathbb{S}^n} f_p\right)^{-\frac{n-4}{4}} - (f(P))^{-\frac{n-4}{4}}\right)$$

$$= O\left(\|f_p - f(P)\|_\infty\right) \quad.$$

On the other hand, by Theorem 2.1, (1.0.3), (2.0.1) and our renormalization

$$\left(\frac{2}{n-4}\right)^{\frac{n}{4}-1} \mathcal{M}_p^{\frac{n}{4}} = E_1[u_p] \geq d_n \left(\fint_{\mathbb{S}^n} |u_p|^{2^\#}\right)^{\frac{2}{2^\#}} = d_n \left(f(P)\right)^{-\frac{2}{2^\#}} \left(\fint_{\mathbb{S}^n} f(P) |u_p|^{2^\#}\right)^{\frac{2}{2^\#}}$$

$$= d_n \left(f(P)\right)^{-\frac{2}{2^\#}} \left(\fint_{\mathbb{S}^n} (f(P) - f_p) |u_p|^{2^\#} + \fint_{\mathbb{S}^n} f_p |u_p|^{2^\#}\right)^{\frac{2}{2^\#}}$$

$$= d_n \left(f(P)\right)^{-\frac{2}{2^\#}} \left(O\left(\|f_p - f(P)\|_\infty\right) + \left(\frac{2}{n-4}\right)^{\frac{n}{4}} \mathcal{M}_p^{\frac{n}{4}}\right)^{\frac{2}{2^\#}} .$$

It follows that $\mathcal{M}_p \geq d_n \left(f(P)\right)^{-\frac{2}{2^\#}} + O\left(\|f_p - f(P)\|_\infty\right)$, and thus we obtain

$$\fint_{\mathbb{S}^n} u_p^2 - \left(\frac{n(n^2-4)}{8f(P)}\right)^{\frac{n-4}{4}} = O\left(\|f_p - f(P)\|_\infty\right) .$$

This concludes the proof of Lemma 5.2. $\square$

**Lemma 5.3.** *There holds*

$$\left\|u_p - \left(\frac{n(n^2-4)}{8f(P)}\right)^{\frac{n-4}{8}}\right\|_{H_2^2}^2 = O\left(\|f_p - f(P)\|_\infty\right) .$$

**Proof.** According to Lemma 5.2, it remains to prove that $\fint_{\mathbb{S}^n} \left|u_p - \left(\frac{n(n^2-4)}{8f(P)}\right)^{\frac{n-4}{8}}\right|^2 = O\left(\|f_p - f(P)\|_\infty\right)$.
We already know from Lemma 5.2 that $\fint_{\mathbb{S}^n} u_p^2 - \left(\frac{n(n^2-4)}{8f(P)}\right)^{\frac{n-4}{4}} = O\left(\|f_p - f(P)\|_\infty\right)$; using the Sobolev-Poincaré inequality, we have

$$\left| \|u_p\|_2 - \left\|\fint_{\mathbb{S}^n} u_p\right\|_2 \right| \leq \left\|u_p - \fint_{\mathbb{S}^n} u_p\right\|_2 \leq C \left(\fint_{\mathbb{S}^n} |\nabla u_p|^2\right)^{\frac{1}{2}} ,$$

and since $\fint_{\mathbb{S}^n} u_p \geq 0$

$$-C \left(\fint_{\mathbb{S}^n} |\nabla u_p|^2\right)^{\frac{1}{2}} + \left(\fint_{\mathbb{S}^n} u_p^2\right)^{\frac{1}{2}} \leq \fint_{\mathbb{S}^n} u_p \leq C \left(\fint_{\mathbb{S}^n} |\nabla u_p|^2\right)^{\frac{1}{2}} + \left(\fint_{\mathbb{S}^n} u_p^2\right)^{\frac{1}{2}} .$$

Hence,

$$-C \|f_p - f(P)\|_\infty^{\frac{1}{2}} \leq \fint_{\mathbb{S}^n} u_p - \left(\frac{n(n^2-4)}{8f(P)}\right)^{\frac{n-4}{8}} \leq C \|f_p - f(P)\|_\infty^{\frac{1}{2}} .$$

Using once again the Sobolev-Poincaré inequality, it is

$$\left(\fint_{\mathbb{S}^n} \left|u_p - \left(\frac{n(n^2-4)}{8f(P)}\right)^{\frac{n-4}{8}} - \left(\fint_{\mathbb{S}^n} u_p - \left(\frac{n(n^2-4)}{8f(P)}\right)^{\frac{n-4}{8}}\right)\right|^2\right)^{\frac{1}{2}} \leq C \left(\fint_{\mathbb{S}^n} |\nabla u_p|^2\right)^{\frac{1}{2}} ,$$

and this implies

$$\left(\fint_{\mathbb{S}^n} \left|u_p - \left(\frac{n(n^2-4)}{8f(P)}\right)^{\frac{n-4}{8}}\right|^2\right)^{\frac{1}{2}} \leq C \left(\fint_{\mathbb{S}^n} |\nabla u_p|^2\right)^{\frac{1}{2}} + \left|\fint_{\mathbb{S}^n} u_p - \left(\frac{n(n^2-4)}{8f(P)}\right)^{\frac{n-4}{8}}\right| .$$

It follows that
$$\fint_{\mathbb{S}^n} \left| u_p - \left(\frac{n(n^2-4)}{8f(P)}\right)^{\frac{n-4}{8}} \right|^2 = O\left(\|f_p - f(P)\|_\infty\right) \quad .$$

This concludes the proof of Lemma 5.3. □

**Lemma 5.4.** $(\Lambda_p)_p$ is bounded.

**Proof.** We have, multiplying (5.1.1) by $u_p \xi_p$ and integrating,
$$\Lambda_p \fint_{\mathbb{S}^n} |u_p|^{2^\#} \xi_p^2 = \frac{n-4}{2} \fint_{\mathbb{S}^n} f_p |u_p|^{2^\#} \xi_p - \fint_{\mathbb{S}^n} (\Delta_h u_p)^2 \xi_p - c_n \fint_{\mathbb{S}^n} |\nabla u_p|^2 \xi_p - \left(\frac{n^2}{2} - c_n \frac{n}{2} + d_n\right) \fint_{\mathbb{S}^n} u_p^2 \xi_p \quad ,$$
and $(u_p)$ is bounded in $H_2^2(\mathbb{S}^n)$; so, as in the proof of Lemma 4.6, we obtain the Lemma. □

**Lemma 5.5.** $(u_p)_p$ is bounded in $L^\infty(\mathbb{S}^n)$.

**Proof.** Consider for $p = \frac{t-1}{t} P \in B$
$$\Omega_p = \left\{ x \in \mathbb{S}^n \mid \left| u_p(x) - \left(\frac{n(n^2-4)}{8f(P)}\right)^{\frac{n-4}{8}} \right| < 1 \right\} \quad .$$

According to lemma 5.3, $|\mathbb{S}^n \setminus \Omega_p| \to 0$ uniformly in $p$ as $\varepsilon_f \to 0$. Now, we write
$$\left(\Delta_h + \frac{c_n}{2}\right)^2 u_p = b_p u_p + w_p$$

where
$$b_p = \left(\frac{n-4}{2} f_p - \Lambda_p \xi_p\right) \left(|u_p|^{2^\# - 2} - \left(\frac{n(n^2-4)}{8f(P)}\right)\right)$$
$$+ \left(\frac{c_n^2}{4} - d_n + \left(\frac{n-4}{2} f_p - \Lambda_p \xi_p\right)\left(\frac{n(n^2-4)}{8f(P)}\right)\right) 1_{\mathbb{S}^n \setminus \Omega_p}$$
$$w_p = \left(\frac{c_n^2}{4} - d_n + \left(\frac{n-4}{2} f_p - \Lambda_p \xi_p\right)\left(\frac{n(n^2-4)}{8f(P)}\right)\right) u_p 1_{\Omega_p} \quad .$$

We have using lemmas 5.3 and 5.4
$$\|b_p\|_{\frac{n}{4}} \leq C \left( \left\| u_p - \frac{n(n^2-4)}{8f(P)} \right\|_{2^\#} + |\mathbb{S}^n \setminus \Omega_p|^{\frac{4}{n}} \right)$$
$$\leq C \left( \|f_p - f(P)\|_\infty^{\frac{1}{2}} + \|f_p - f(P)\|_\infty^{\frac{2}{n}} \right) \quad ,$$

where $C$ is independent of $p$. It follows, using Remark 3.2 and Proposition 3.1, since $b_p \in L^\infty(\mathbb{S}^n)$ and $|\mathbb{S}^n \setminus \Omega_p| \to 0$ as $\varepsilon_f \to 0$ uniformly in $p$, that for $\varepsilon_f$ small enough, we have for all $p \in B$
$$\|u_p\|_\infty \leq C \|w_p\|_\infty \quad ,$$
where $C$ is independent of $p$. Thanks to Lemma 5.4 $(w_p)$ is bounded in $L^\infty(\mathbb{S}^n)$ uniformly in $p$. This concludes the proof of the Lemma. □

**Remarks 5.6.** From the previous Lemma, the proof of Lemma 5.2, and the proof of Lemma 5.3, we can deduce the following estimate
$$\left\| u_p - \left(\frac{n(n^2-4)}{8f(P)}\right)^{\frac{n-4}{8}} \right\|_{H_2^2}^2 = O\left(\|f_p - f(P)\|_1\right) \quad .$$

**Lemma 5.7.** *There holds*
$$\Lambda_p = O\left(\|f_p - f(P)\|_\infty\right) \quad.$$

**Proof.** From Theorem 1.1, we deduce
$$\Lambda_p \fint_{\mathbb{S}^n} |\nabla \xi_p|^2 |u_p|^{2^\#} = \frac{n-4}{8} \fint_{\mathbb{S}^n} <\nabla f_p, \nabla \xi_p> |u_p|^{2^\#} \quad.$$

Writing
$$\fint_{\mathbb{S}^n} <\nabla f_p, \nabla \xi_p> |u_p|^{2^\#} = \fint_{\mathbb{S}^n} <\nabla (f_p - f(P)), \nabla \xi_p> |u_p|^{2^\#}$$
$$= n \fint_{\mathbb{S}^n} (f_p - f(P)) \xi_p |u_p|^{2^\#} - 2^\# \fint_{\mathbb{S}^n} (f_p - f(P)) <\nabla u_p, \nabla \xi_p> |u_p|^{2^\# - 2} u_p \quad,$$

and using Hölder inequality, Lemma 5.3 and Lemma 5.5, we get
$$\Lambda_p \fint_{\mathbb{S}^n} |\nabla \xi_p|^2 |u_p|^{2^\#} = O\left(\|f_p - f(P)\|_\infty\right) \quad.$$

So, as in the proof of Lemma 4.6, we derive Lemma 5.7. $\square$

**Lemma 5.8.** *There holds*
$$\left\| u_p - \left(\frac{n(n^2 - 4)}{8 f(P)}\right)^{\frac{n-4}{8}} \right\|_{C^2} = O\left(\|f_p - f(P)\|_\infty\right) \quad.$$

**Proof.** Consider for $C$ large, and for $p \in B$
$$\Omega_p = \left\{ x \in \mathbb{S}^n \mid \left| u_p(x) - \left(\frac{n(n^2 - 4)}{8 f(P)}\right)^{\frac{n-4}{8}} \right| < C \|f_p - f(P)\|_\infty \right\} \quad.$$

According to lemma 5.3, $|\mathbb{S}^n \setminus \Omega_p| \to 0$ as $C \to +\infty$ uniformly in $p$. We have

(5.8.1) $$\left(\Delta_h + \frac{c_n}{2}\right)^2 \left( u_p - \left(\frac{n(n^2 - 4)}{8 f(P)}\right)^{\frac{n-4}{8}} \right) = b_p \left( u_p - \left(\frac{n(n^2 - 4)}{8 f(P)}\right)^{\frac{n-4}{8}} \right) + w_p \quad,$$

where
$$b_p = \left(\frac{n-4}{2} f_p |u_p|^{2^\# - 2} - d_n + \frac{c_n^2}{4}\right) 1_{\mathbb{S}^n \setminus \Omega_p}$$
$$w_p = \left(\frac{n-4}{2} f_p |u_p|^{2^\# - 2} - d_n + \frac{c_n^2}{4}\right) \left( u_p - \left(\frac{n(n^2 - 4)}{8 f(P)}\right)^{\frac{n-4}{8}} \right) 1_{\Omega_p} - \Lambda_p \xi_p |u_p|^{2^\# - 2} u_p$$
$$+ \left(\frac{n-4}{2} f_p |u_p|^{2^\# - 2} - d_n\right) \left(\frac{n(n^2 - 4)}{8 f(P)}\right)^{\frac{n-4}{8}} \quad,$$

we clearly get, using Lemma 5.5
$$\|b_p\|_{\frac{n}{4}} \leq |\mathbb{S}^n \setminus \Omega_p|^{\frac{4}{n}} \quad.$$

It follows that for $C$ large enough (independent of $p$), $\|b_p\|_{\frac{n}{4}} \leq \delta^\infty$ for all $p \in B$ (here $\delta^\infty$ is given by Remark 3.2). From this we deduce, using Remark 3.3, Proposition 3.1, Lemma 5.5, Lemma 5.7 and (5.8.1), that for all $p \in B$

$$(5.8.2) \qquad \left\| u_p - \left(\frac{n(n^2-4)}{8f(P)}\right)^{\frac{n-4}{8}} \right\|_\infty \leq C \|w_p\|_{\frac{n}{4}} \leq C_1 \|f_p - f(P)\|_\infty \quad,$$

where $C_1$ is independent of $p$. Then using (5.8.1) and (5.8.2), it easily follows that for all $p \in B$

$$\left\| u_p - \left(\frac{n(n^2-4)}{8f(P)}\right)^{\frac{n-4}{8}} \right\|_{C^2} \leq C \|f_p - f(P)\|_\infty \quad,$$

where $C$ is independent of $p$. This concludes the proof of the Lemma. $\square$

**Remarks 5.9.** It follows from the previous Lemma that if $\|f_p - f(P)\|_\infty$ is small enough (this smallness being independent of $t$ and $P$), $u_p > 0$ on $\mathbb{S}^n$. So, if for some $t_0$ it is $\Lambda_{t_0} = 0$, this gives rise, according to our transformation rules, to a positive solution of the original equation $(E)$.

## §6. Continuous dependence on the parameter.

In this section we prove that, given $p \in B$, $p = (P, t)$, the function $u_p = u_{P,t}$ given by Lemma 5.1 is uniquely determined and that $u_p$, as well as the Lagrange multiplier $\Lambda_p$, vary continuously in $B$. For $p = \frac{t-1}{t}P$, $P \in \mathbb{S}^n$, we consider the functional

$$\bar{J}_p[u] = \frac{\fint_{\mathbb{S}^n} (\Delta_h u)^2 + c_n \fint_{\mathbb{S}^n} |\nabla u|^2 + d_n \fint_{\mathbb{S}^n} u^2}{\left(\fint_{\mathbb{S}^n} f_p |u|^{2^\#}\right)^{\frac{2}{2^\#}}} \quad,$$

where $f_p$ stands for $f_p$, and we let

$$\mathcal{M}_p = \inf_{u \in \mathcal{S}_{2^\#}} \bar{J}_p[u] \quad.$$

**Proposition 6.1.** For $\varepsilon_f = \left\|f - \frac{n(n^2-4)}{8}\right\|_\infty$ sufficiently small, the functional $\bar{J}_p$ has a unique minimum $u_p$ in the class $\mathcal{S}_{2^\#}^0$. The map $p \to u_p$ is continuous from $B$ into $\mathcal{S}_{2^\#}^0$.

**Proof.** To verify the uniqueness assertion, we assume that there exists $p \in B$ for which $\bar{J}_p$ has two distinct minima in $\mathcal{S}_{2^\#}^0$, $u_0$ and $u_1$. According to remark 5.9 we can choose $\varepsilon_f$ small enough such that $u_1$ and $u_2$ are positive. For convenience reason, by a simple renormalization, (and without loss of generality) we assume that $u_0$ and $u_1$ are solutions respectively of

$$\Delta_h^2 u_0 + c_n \Delta_h u_0 + d_n u_0 = \frac{n-4}{2} f_p u_0^{2^\#-1} - \Lambda_0 \xi_0 u_0^{2^\#-1}$$

and

$$\Delta_h^2 u_1 + c_n \Delta_h u_1 + d_n u_1 = \frac{n-4}{2} f_p u_1^{2^\#-1} - \Lambda_1 \xi_1 u_1^{2^\#-1} \quad,$$

where $\Lambda_0, \Lambda_1 \in \mathbb{R}^+$, $\xi_0, \xi_1 \in \mathbb{H}$, $\|\xi_0\|_\infty = \|\xi_1\|_\infty = 1$. We set, for $\lambda \in [0; 1]$,

$$u_\lambda^{2^\#} = \lambda u_0^{2^\#} + (1-\lambda) u_1^{2^\#} \quad.$$

For each $\lambda$ there holds

$$\begin{cases} \dot{u}_\lambda = \frac{1}{2^\#} u_\lambda^{1-2^\#} \left(u_0^{2^\#} - u_1^{2^\#}\right) \quad; \\ \ddot{u}_\lambda = -\left(2^\# - 1\right) u_\lambda^{-1} (\dot{u}_\lambda)^2 \quad. \end{cases}$$

Also, differentiating with respect to $\lambda$, we deduce

$$u_\lambda \in \mathcal{S}_{2^\#} \Rightarrow \fint_{\mathbb{S}^n} u_\lambda^{2^\#-1} \dot{u}_\lambda \xi = 0 \qquad \text{for all } \xi \in \mathbb{H} \quad .$$

Hence, by Hölder's inequality and Lemma 5.8, for all $\xi \in \mathbb{H}$

$$\left(\frac{n(n^2-4)}{8f(P)}\right)^{\frac{n+4}{8}} \left|\fint_{\mathbb{S}^n} \dot{u}_\lambda \xi\right| = \left|\fint_{\mathbb{S}^n} \left(\left(\frac{n(n^2-4)}{8f(P)}\right)^{\frac{n+4}{8}} - u_\lambda^{2^\#-1}\right) \dot{u}_\lambda \xi\right| = O(\varepsilon_f) \|\dot{u}_\lambda\|_2 \|\xi\|_2$$

Decomposing $\dot{u}_\lambda$ into $\alpha_\lambda + \Psi_\lambda + \xi_\lambda$, where $\alpha_\lambda \in \mathbb{R}$, $\xi_\lambda$ is the orthogonal projection of $\dot{u}_\lambda$ onto $\mathbb{H}$ and $\Psi_\lambda$ is orthogonal to $\mathbb{R} \oplus \mathbb{H}$, we have

$$\fint_{\mathbb{S}^n} |\nabla \dot{u}_\lambda|^2 = \fint_{\mathbb{S}^n} |\nabla \Psi_\lambda|^2 + \fint_{\mathbb{S}^n} |\nabla \xi_\lambda|^2 + 2\fint_{\mathbb{S}^n} \nabla \Psi_\lambda \nabla \xi_\lambda = \fint_{\mathbb{S}^n} |\nabla \Psi_\lambda|^2 + n\fint_{\mathbb{S}^n} \xi_\lambda^2 \quad .$$

On the other hand

$$\fint_{\mathbb{S}^n} |\nabla \dot{u}_\lambda|^2 = \fint_{\mathbb{S}^n} \dot{u}_\lambda (\Delta \Psi_\lambda + n\xi_\lambda) = \fint_{\mathbb{S}^n} |\nabla \Psi_\lambda|^2 + n\fint_{\mathbb{S}^n} \xi_\lambda \dot{u}_\lambda \quad .$$

Taking into account that $\fint_{\mathbb{S}^n} \dot{u}_\lambda \xi_\lambda = O(\varepsilon_f) \|\dot{u}_\lambda\|_2 \|\xi_\lambda\|_2$, it is

$$\fint_{\mathbb{S}^n} \xi_\lambda^2 = O(\varepsilon_f) \|\dot{u}_\lambda\|_2^2 \quad .$$

It follows, from the Courant-Fischer characterization of the eigenvalues and from Theorem 2.2, that

$$(6.1.1) \quad \begin{aligned} \fint_{\mathbb{S}^n} (\Delta_h \dot{u}_\lambda)^2 + c_n \fint_{\mathbb{S}^n} |\nabla \dot{u}_\lambda|^2 &\geq \left(4(n+1)^2 + c_n 2(n+1)\right) \fint_{\mathbb{S}^n} \Psi_\lambda^2 + \varepsilon_f O(\|\dot{u}_\lambda\|_2) \\ &= \left(4(n+1)^2 + c_n 2(n+1) + O(\varepsilon_f)\right) \fint_{\mathbb{S}^n} \dot{u}_\lambda^2 \quad . \end{aligned}$$

Moreover, setting $f_0 = \frac{n(n^2-4)}{8}$, we have the following estimates

$$(i) \quad \fint_{\mathbb{S}^n} f_p u_\lambda^{2^\#} = \fint_{\mathbb{S}^n} (f_p - f_0) u_\lambda^{2^\#} + \fint_{\mathbb{S}^n} f_0 u_\lambda^{2^\#} = f_0 (1 + O(\varepsilon_f)) \qquad \text{by Lemma 5.3} \quad ;$$

From $\fint_{\mathbb{S}^n} u_\lambda^{2^\#} = constant$ we deduce $\fint_{\mathbb{S}^n} u_\lambda^{2^\#-1} \dot{u}_\lambda = 0$ and hence

$$(ii) \quad \fint_{\mathbb{S}^n} f_p u_\lambda^{2^\#-1} \dot{u}_\lambda = O(\varepsilon_f) \left(\fint_{\mathbb{S}^n} \dot{u}_\lambda^2\right)^{\frac{1}{2}} \quad ;$$

and using Lemma 5.3 and Lemma 5.8

$$(iii) \quad \fint_{\mathbb{S}^n} (\Delta_h u_\lambda)^2 + c_n \fint_{\mathbb{S}^n} |\nabla u_\lambda|^2 + d_n \fint_{\mathbb{S}^n} u_\lambda^2 = \frac{n-4}{2} f_0 (1 + O(\varepsilon_f)) \quad ;$$

$$(iv) \quad \fint_{\mathbb{S}^n} f_p u_\lambda^{2^\#-2} \dot{u}_\lambda^2 = \fint_{\mathbb{S}^n} (f_p - f_0) u_\lambda^{2^\#-2} \dot{u}_\lambda^2 + \fint_{\mathbb{S}^n} f_0 u_\lambda^{2^\#-2} \dot{u}_\lambda^2 = f_0 (1 + O(\varepsilon_f)) \left(\fint_{\mathbb{S}^n} \dot{u}_\lambda^2\right) \quad .$$

We have also

$$d^2 \bar{J}_p[u].(\varphi, \Psi)$$
$$= 2 \left( \fint_{\mathbb{S}^n} f_p u^{2^\#} \right)^{-\frac{2}{2^\#}-1} \times \left\{ \fint_{\mathbb{S}^n} \left( f_p u^{2^\#} \right) . \left( \fint_{\mathbb{S}^n} \Delta_h \varphi \Delta_h \Psi + c_n \fint_{\mathbb{S}^n} \nabla \varphi \nabla \Psi + d_n \fint_{\mathbb{S}^n} \varphi \Psi \right) \right.$$
$$-2 \left( \fint_{\mathbb{S}^n} f_p u^{2^\#-1} \Psi \right) \left( \fint_{\mathbb{S}^n} \Delta u \Delta \varphi + c_n \fint_{\mathbb{S}^n} \nabla u \nabla \varphi + d_n \fint_{\mathbb{S}^n} u \varphi \right)$$
$$-2 \left( \fint_{\mathbb{S}^n} f_p u^{2^\#-1} \varphi \right) \left( \fint_{\mathbb{S}^n} \Delta u \Delta \Psi + c_n \fint_{\mathbb{S}^n} \nabla u \nabla \Psi + d_n \fint_{\mathbb{S}^n} u \Psi \right)$$
$$-(2^\# - 1) \left( \fint_{\mathbb{S}^n} (\Delta_h u)^2 + c_n \fint_{\mathbb{S}^n} |\nabla u|^2 + d_n \fint_{\mathbb{S}^n} u^2 \right) \left( \fint_{\mathbb{S}^n} f_p u^{2^\#-2} \varphi \Psi \right)$$
$$+ \left( \frac{2}{2^\#} + 1 \right) \left( \fint_{\mathbb{S}^n} (\Delta_h u)^2 + c_n \fint_{\mathbb{S}^n} |\nabla u|^2 + d_n \fint_{\mathbb{S}^n} u^2 \right) \times$$
$$\left. \left( \fint_{\mathbb{S}^n} f_p u^{2^\#-1} \varphi \right) \left( \fint_{\mathbb{S}^n} f_p u^{2^\#-1} \Psi \right) \left( \fint_{\mathbb{S}^n} f_p u^{2^\#} \right)^{-1} \right\} \quad .$$

It is clear, using the Sobolev inequality, that

$$u \to D^2 \bar{J}_p[u](.,.) \quad \text{is a continuous map from } H_2^2(\mathbb{S}^n) \text{ to } L\left( H_2^2(\mathbb{S}^n) \times H_2^2(\mathbb{S}^n), \mathbb{R} \right) \quad .$$

We have

$$\frac{d^2}{d\lambda^2} \bar{J}_p[u_\lambda] = d^2 \bar{J}_p[u_\lambda].(\dot{u}_\lambda, \dot{u}_\lambda) + d\bar{J}_p[u_\lambda].\ddot{u}_\lambda \quad ,$$

and moreover from $(i) - (iv)$ we have

$$d^2 \bar{J}_p[u_\lambda].(\dot{u}_\lambda, \dot{u}_\lambda)$$
$$= 2 \left( \fint_{\mathbb{S}^n} f_p u_\lambda^{2^\#} \right)^{-\frac{2}{2^\#}-1} \times \left\{ (1 + O(\varepsilon_f)) f_0 \fint_{\mathbb{S}^n} (\Delta_h \dot{u}_\lambda)^2 + (1 + O(\varepsilon_f)) f_0 c_n \fint_{\mathbb{S}^n} |\nabla \dot{u}_\lambda|^2 \right.$$
$$\left. + d_n f_0 \left( 1 - (2^\# - 1) \right) \fint_{\mathbb{S}^n} \dot{u}_\lambda^2 (1 + O(\varepsilon_f)) \right\}$$
$$= 2 \left( \fint_{\mathbb{S}^n} f_p u_\lambda^{2^\#} \right)^{-\frac{2}{2^\#}-1} f_0 (1 + O(\varepsilon_f)) \times \left\{ \fint_{\mathbb{S}^n} (\Delta_h \dot{u})^2 + c_n \fint_{\mathbb{S}^n} |\nabla \dot{u}_\lambda|^2 - \frac{8}{n-4} d_n \fint_{\mathbb{S}^n} \dot{u}_\lambda^2 \right\} \quad ,$$

and

$$d\bar{J}_p[u_\lambda].\ddot{u}_\lambda = 2 \left( \fint_{\mathbb{S}^n} f_p u_\lambda^{2^\#} \right)^{-\frac{2}{2^\#}} \times \left\{ \fint_{\mathbb{S}^n} \Delta_h u_\lambda \Delta_h \ddot{u}_\lambda + c_n \fint_{\mathbb{S}^n} \nabla u_\lambda \nabla \ddot{u}_\lambda + d_n \fint_{\mathbb{S}^n} u_\lambda \ddot{u}_\lambda \right.$$
$$\left. - \left( \fint_{\mathbb{S}^n} (\Delta_h u_\lambda)^2 + c_n \fint_{\mathbb{S}^n} |\nabla u_\lambda|^2 + d_n \fint_{\mathbb{S}^n} u_\lambda^2 \right) \fint_{\mathbb{S}^n} f_p u_\lambda^{2^\#-1} \ddot{u}_\lambda \left( \fint_{\mathbb{S}^n} f_p u_\lambda^{2^\#} \right)^{-1} \right\} \quad .$$

One can check with some straightforward computations and using Lemma 5.8, that

$$d\bar{J}_p[u_\lambda].\ddot{u}_\lambda = 2 \left( \fint_{\mathbb{S}^n} f_p u_\lambda^{2^\#} \right)^{-\frac{2}{2^\#}} \left\{ (1 - 2^\#) \left( \fint_{\mathbb{S}^n} \Delta_h u_\lambda \Delta_h \left( \frac{\dot{u}_\lambda^2}{u_\lambda} \right) + c_n \fint_{\mathbb{S}^n} \nabla u_\lambda \nabla \left( \frac{\dot{u}_\lambda^2}{u_\lambda} \right) + d_n \fint_{\mathbb{S}^n} u_\lambda \frac{\dot{u}_\lambda^2}{u_\lambda} \right) \right.$$
$$\left. - \left( \fint_{\mathbb{S}^n} (\Delta_h u_\lambda)^2 + c_n \fint_{\mathbb{S}^n} |\nabla u_\lambda|^2 + d_n \fint_{\mathbb{S}^n} u_\lambda^2 \right) \fint_{\mathbb{S}^n} f_p u_\lambda^{2^\#-1} \ddot{u}_\lambda \left( \fint_{\mathbb{S}^n} f_p u_\lambda^{2^\#} \right)^{-1} \right\}$$
$$= O(\varepsilon_f) \left( \fint_{\mathbb{S}^n} (\Delta_h \dot{u}_\lambda)^2 + \fint_{\mathbb{S}^n} |\nabla \dot{u}_\lambda|^2 + \fint_{\mathbb{S}^n} \dot{u}_\lambda^2 \right) \quad .$$

It follows from (6.1.1) that

$$\frac{d^2}{d\lambda^2} \bar{J}_p[u_\lambda] \geq C \left\{ \left(1 - \left(\frac{n^3 - 4n}{2(n^3 + 3n^2 + 2n)}\right) + O(\varepsilon_f)\right) \left(\fint_{\mathbb{S}^n} (\Delta_h \dot{u}_\lambda)^2 + c_n \fint_{\mathbb{S}^n} |\nabla \dot{u}_\lambda|^2\right) \right\} > 0$$

for $\varepsilon_f$ small. This means that $\bar{J}_p[u_\lambda]$ is a strictly convex function on $[0;1]$, hence this contradicts the assumption that $u_0$ and $u_1$ are minima unless $u_0 = u_1$.

Let us observe that, reasoning as above, one can prove that $d^2 \bar{J}_p[u_p]$ is positive definite on the tangent space of the constraint $\mathcal{S}_{2\#}$, so using the Implicit Functions Theorem, it follows that $p \to u_p$ is continuous This concludes the proof of Proposition 6.1. □

## §7. Comparison of the maps $\Lambda$ and $G$.

Consider the map $\Lambda : B \to \mathbb{R}$. We can associate to this map the folowing map $\overrightarrow{\Lambda} : B \to \mathbb{R}^{n+1}$ consisting in the components of $\Lambda$ in the basis $\overrightarrow{\xi}$. Namely, in the basis $\overrightarrow{\xi}$ of $\mathbb{H}$ we can write

$$\Lambda = \sum_{i=1}^{n+1} \Lambda_i \xi_i \quad,$$

and then

$$\overrightarrow{\Lambda} = (\Lambda_1, \ldots, \Lambda_{n+1}) \quad.$$

In this section we prove that $\overrightarrow{\Lambda}$ has the same degree as the map $G$ defined in the introduction. We set

$$\begin{cases} C_{i,j}(p) = \fint_{\mathbb{S}^n} <\nabla \xi_i, \nabla \xi_j> u_p^{2^\#} \\ A_j(p) = \fint_{\mathbb{S}^n} <\nabla f, \nabla \xi_j> u_p^{2^\#} \end{cases}$$

We observe that, since $u_p$ is close to 1 by Lemma 5.8, for $\varepsilon_f$ small, $(C_{i,j}(p))$ is a positive definite matrix. Therefore the Kazdan-Warner condition can be rewritten as

$$\overrightarrow{\Lambda}(p) = C(p)^{-1} A(p)$$

It follows immediately that $\overrightarrow{\Lambda}$ is continuous (using Proposition 6.1) and

**Lemma 7.1.**

(i) $\overrightarrow{\Lambda}(p) = 0 \Leftrightarrow A(p) = 0 \quad;$

(ii) $deg(\overrightarrow{\Lambda}|\{(P,t), t < t_0\}, 0) = deg(A|\{(P,t), t < t_0\}, 0) \quad \text{for all } t_0 > 1 \quad.$

Now consider (using the notation of the introduction)

$$A(P,t) = \fint_{\mathbb{S}^n} <\nabla (f \circ \varphi_{P,t}), \nabla \overrightarrow{\xi}> u_p^{2^\#} \quad,$$

$$G_1(P,t) = n \left(\frac{n(n^2-4)}{8f(P)}\right)^{\frac{n}{4}} \fint_{\mathbb{S}^n} (f \circ \varphi_{P,t}) \overrightarrow{\xi} = n \left(\frac{n(n^2-4)}{8f(P)}\right)^{\frac{n}{4}} G(P,t) \quad.$$

(this means that the components of $A(P,t)$ are $\fint_{\mathbb{S}^n} <\nabla (f \circ \varphi_{P,t}), \nabla \xi_i> u_p^{2^\#}$, for $i = 1, \ldots, n+1$). We choose $\varepsilon$ sufficiently small so that the error term in $f(y) = f(P) + \sum_{k=1}^{\alpha} f_k(y) + O\left(|y|^{\alpha+1}\right)$ is small, say $O\left(|y|^{\alpha+1}\right) \leq C |y|^{\alpha+\frac{1}{2}}$ for $|y| < \varepsilon$. We have

$$A(P,t) = G_1(P,t) + I + II \quad,$$

where

$$I = n \fint_{\mathbb{S}^n} (f \circ \varphi_{P,t} - f(P)) \vec{\xi} \left( u_p^{2^\#} - \left( \frac{n(n^2-4)}{8f(P)} \right)^{\frac{n}{4}} \right) \quad ;$$

$$II = \fint_{\mathbb{S}^n} (f \circ \varphi_{P,t} - f(P)) <\nabla \vec{\xi}, \nabla (u_p^{2^\#})> \quad .$$

We set, using stereographic coordinates, $\Omega_t = \{y \in \mathbb{R}^n | \, |y| \leq \varepsilon t\}$ and we have $|\mathbb{S}^n \setminus \Omega_t| = O\left(\frac{1}{t^n}\right)$. Clearly we have

$$\left| \fint_{\Omega_t} (f \circ \varphi_{P,t} - f(P)) \vec{\xi} \left( u_p^{2^\#} - \left( \frac{n(n^2-4)}{8f(P)} \right)^{\frac{n}{4}} \right) \right| \leq \quad ,$$

$$\left( \fint_{\Omega_t} (f \circ \varphi_{P,t} - f(P))^2 \right)^{\frac{1}{2}} \left( \fint_{\Omega_t} \left( u_p^{2^\#} - \left( \frac{n(n^2-4)}{8f(P)} \right)^{\frac{n}{4}} \right)^2 \right)^{\frac{1}{2}} \quad .$$

and

$$\left| \fint_{\Omega_t} (f \circ \varphi_{P,t} - f(P)) <\nabla x, \nabla (u_p^{2^\#})> \right| \leq$$

$$\left( \fint_{\Omega_t} (f \circ \varphi_{P,t} - f(P))^2 \right)^{\frac{1}{2}} \left( \fint_{\Omega_t} \left| \nabla (u_p^{2^\#}) \right|^2 \right)^{\frac{1}{2}} \quad .$$

To have further control on the error terms we point out the following estimates.

**Lemma 7.2.** *If $P$ is a critical point of $f$ non degenerated of order $\alpha$, we have the following estimates*

$$\fint_{\mathbb{S}^n} (f_p - f(P))^2 = \begin{cases} O\left(\frac{1}{t^{2\alpha}}\right) & \text{if} \quad 2\alpha < n \\ O\left(\frac{1}{t^n} \log t\right) & \text{if} \quad 2\alpha = n \\ O\left(\frac{1}{t^n}\right) & \text{if} \quad 2\alpha > n \end{cases}$$

$$\fint_{\mathbb{S}^n} |\nabla u_p|^2 = O\left(\|f_p - f(P)\|_1\right) = \begin{cases} O\left(\frac{1}{t^\alpha}\right) & \text{if} \quad \alpha < n \\ O\left(\frac{1}{t^n} \log t\right) & \text{if} \quad \alpha = n \end{cases}$$

$$\left\| u_p^{2^\#} - \left( \frac{n(n^2-4)}{8f(P)} \right)^{\frac{n}{4}} \right\|_2^2 = O\left(\|f_p - f(P)\|_1\right) = \begin{cases} O\left(\frac{1}{t^\alpha}\right) & \text{if} \quad \alpha < n \\ O\left(\frac{1}{t^n} \log t\right) & \text{if} \quad \alpha = n \end{cases}$$

The proof of the first estimate of Lemma 7.2 is obtained by a simple computation, the second and the third estimates are obtained using remark 5.6 and a simple computation.

**Proposition 7.3.** *Assume that $f$ is everywhere positive, uniformly non degenerated of order at most $n$ if $n$ is even, and of order at most $n-1$ if $n$ is odd, . Then we have for all $P \in \mathbb{S}^n$ and for all $t$ large (with the notation of this section)*

$$G(P,t).A(P,t) > 0 \quad .$$

Proposition 7.3 is a direct consequence of Lemma 7.2.

One has to notice that in the case $n = 2k+1 = \alpha$, Lemma 7.2 shows that the corrections terms $I$ and $II$ may in fact be the dominant one. For this reason we exclude this case from consideration and this accounts for the difference between the case where the dimension of the manifold is odd and the case where this dimension is even in the statement of Theorem 1.3.

## §8. Proof of Theorem 1.3 and of Corollary 1.4.

By a simple compactness argument and using Proposition 7.3, there exists $t_0$ such that for all $t \geq t_0$

$$G(P,t).A(P,t) > 0 \quad .$$

Then, for $0 \leq r \leq 1$,

$$deg\left(rG(P,t) + (1-r)A(P,t), \{p \mid t < t_0\}, 0\right) = constant \quad .$$

It follows that, using Lemma 7.1 and under the notation of section 7

$$deg\left(\overrightarrow{\Lambda}(p), \{p \mid t < t_0\}, 0\right) = deg\left(A(p), \{p \mid t < t_0\}, 0\right) = deg\left(G(p), \{p \mid t < t_0\}, 0\right) \neq 0 \quad .$$

Hence, there exists $p \in B$ such that $\overrightarrow{\Lambda}(p) = 0$. This concludes the proof of Theorem 1.3. □

In order to show Corollary 1.4, we recall the formulas for the stereographic projection: the point $(x, x_{n+1}) \in \mathbb{S}^n \subset \mathbb{R}^{n+1}$ is projected through the north pole on the point $y \in \mathbb{R}^n$ by the following formulas

(8.0.1) $$x = \frac{2y}{|y|^2 + 1} \quad ; \quad x_{n+1} = \frac{|y|^2 - 1}{|y|^2 + 1} \quad .$$

Consider the stereographic projection $\pi_P : \mathbb{S}^n \to \mathbb{R}^n$ through the point $-P$, and define $\widetilde{f} : \mathbb{R}^n \to \mathbb{R}$ in the following way

$$\widetilde{f}(y) = f\left(\pi_P^{-1}(y)\right) \quad .$$

Then, using formulas (8.0.1) one can prove that

$$G(P,t) = \fint_{\mathbb{S}^n} \widetilde{f}\left(\frac{1}{t}\right)\left(1 + |y|^2\right)^{-(n+1)} \begin{pmatrix} 2y \\ 1 - |y|^2 \end{pmatrix} dy \quad .$$

Here the quantity $\begin{pmatrix} 2y \\ 1 - |y|^2 \end{pmatrix} dy$ is considered as an $(n+1)$-tuple. Hence expanding the above expression in powers of $\frac{1}{t}$, one finds

$$G(P,t) = \begin{pmatrix} a_1 \nabla \widetilde{f}(P)\frac{1}{t} \\ a_2 \Delta \widetilde{f}(P)\frac{1}{t^2} \end{pmatrix} + o\left(\frac{1}{t^2}\right) \quad ,$$

where $a_1$ and $a_2$ are non-zero coefficients given by

$$a_1 = \frac{2}{n}\fint_{\mathbb{S}^n} |y|^2 \left(1 + |y|^2\right)^{-(n+1)} \quad ; \quad a_2 = \frac{2}{n}\fint_{\mathbb{S}^n} \left(|y|^2 - |y|^4\right)\left(1 + |y|^2\right)^{-(n+1)} \quad ,$$

and where $o\left(\frac{1}{t^2}\right)$ is uniform in $P$. So under the assumptions of the Corollary, we deduce that $f$ is uniformly non degenerated of order 2. Moreover, the arguments of Chang-Gursky-Yang [15] show that under the condition

$$\sum_{P \in \mathbb{S}^n \mid \nabla f(P)=0 \text{ and } \Delta_h f(P) > 0} (-1)^{m(P,f)} \neq -1 \quad ,$$

the degree is different from zero. This concludes the proof of Corollary 1.4. □

**Z.D.** : Université de Cergy-Pontoise - Département de Mathématiques - Site de Saint-Martin - 2 Avenue Adolphe Chauvin - F 95302 Cergy-Pontoise Cedex, France

**A.M.** : Rutgers University - Department of Mathematics - Hill Center-Busch Campus - Rutgers, The State University of New Jersey - 110 Frelinghuysen Rd - Piscataway, NJ 08854-8019, USA

**M.O.A.** : Scuola Internazionale Superiore di Studi Avanzati (SISSA) - Via Beirut, 2-4 - 34013 Trieste - Italy

e-mail : Zindine.Djadli@math.u-cergy.fr  ,  malchiod@math.rutgers.edu  ,  ahmedou@sissa.it